\documentclass[12pt]{article}
\usepackage{mathrsfs}
\usepackage{amsfonts}
\usepackage{graphicx}
\begin{document}
\title{The distance between two separating, reducing slopes is at most 4}
\author {Mingxing Zhang, Ruifeng Qiu and Yannan Li\thanks{This work is supported by NSFC}}
\date{}
\maketitle
\begin{abstract}
Let $M$ be a simple 3-manifold such that one component of $\partial
M$, say $F$, has genus at least two. For a slope $\alpha$ on $F$, we
denote by $M(\alpha)$ the manifold obtained by attaching  a 2-handle
to $M$ along a regular neighborhood of $\alpha$ on $F$. If
$M(\alpha)$ is reducible, then $\alpha$ is called a reducing slope.
In this paper, we shall prove that the distance between two
separating, reducing slopes on $F$ is at most 4.


\end{abstract}

\section{Introduction}
\label{intro} Let $M$ be a compact, orientable 3-manifold such that
$\partial M$ contains no spherical components. $M$ is said to be
simple if $M$ is irreducible, $\partial$-irreducible, anannular and
atoroidal.

Let $M$ be a simple 3-manifold. For a component $F$ of $\partial M$,
a slope $\gamma$ on $F$ is an isotopy class of essential simple
closed curves on $F$. For a slope $\gamma$ on $F$, we denote by
$M(\gamma)$ the manifold obtained by attaching a 2-handle to $M$
along a regular neighborhood of $\gamma$ on $F$, then capping off a
possible 2-sphere component of the resulting manifold by a 3-ball. A
slope $\gamma$ on $F$ is said to be reducing if $M(\gamma)$ is
reducible. The distance between two slopes $\alpha$ and $\beta$ on
$F$, denoted by $\Delta(\alpha,\beta)$, is the minimal geometric
intersection number among all the curves representing the slopes.
Note that if $F$ is a torus, then $M(\gamma)$ is the Dehn filling
along $\gamma$. Two important results about reducing handle
additions on simple 3-manifolds are the following:

(1) Suppose that $F$ is a torus, $\alpha$ and $\beta$ are two
reducing slopes on $F$.  Gordon and Luecke[GL1] proved that
$\Delta(\alpha,\beta)\leq 1$. This means that there are at most
three reducing slopes on $F$.

(2) Suppose that $g(F)>1$. Scharlemann and Wu[SW] proved that there
are only finitely many basic degenerating slopes on $F$. As a
corollary of this result, there are only finitely many separating,
reducing slopes on $F$.

In this paper, we shall continue to study reducing handle additions.
The main result is the following theorem:

{\bf Theorem 1.} \ Suppose that $M$ is a simple 3-manifold, and $F$
is a genus at least two component of $\partial M$. If $\alpha$ and
$\beta$ are two separating, reducing slopes on $F$, then
$\Delta(\alpha,\beta)\leq 4$.

{\bf Comments on Theorem 1.}

{\bf 1.} \ It is possible that $\Delta(\alpha,\beta)$ is arbitrarily
large when $\alpha$ and $\beta$ are two non-separating, reducing
slopes on $F$. For example, one can construct a simple 3-manifold
$N$ such that there is a separating, reducing slope $\gamma$ on
$\partial N$ which bounds a punctured torus $T$ in $\partial N$.
Then $N(\gamma)$ is reducible and $\partial N(\gamma)$ contains a
toral component $T^{*}$ such that $T\subset T^{*}$. By the [GL2] and
[SW], there are infinitely many slopes $\alpha$ on $T$ such that
$N(\alpha)=N(\gamma)(\alpha)$ is reducible.

{\bf 2.} \ Let $M$ be a simple 3-manifold containing no essential
closed surfaces of genus $g$. Suppose that $\alpha$ and $\beta$ are
separating slopes on $\partial M$ such that $M(\alpha)$ and
$M(\beta)$ contains an essential closed surface of genus $g$. If
$g\leq1$, then $\Delta(\alpha,\beta)\leq14$, see [SW]. If $g>1$,
then it is possible that $\Delta(\alpha,\beta)$ is arbitrarily
large, see [QW1] and [QW2].

\section{Labeled graph}
\label{labelgraph}

The following Lemma follows from the proof of Lemma 3.3 in [SW].

{\bf Lemma 2.1.} \ Suppose $M$ is a simple manifold. If $\alpha$ is
a separating, reducing slope and $M(\alpha)$ is
$\partial$-irreducible, then $M$ contains an incompressible and
$\partial$-incompressible planar surface in $M$ with all boundary
components having the same slope $\alpha$. \label{lemmaplanar}

{\bf Proof.} \ Suppose $P$ is a planar surface in $M$ with all
boundary components parallel to $\alpha$. Capping off all such
components by mutually disjoint disks in $M(\alpha)$, we get a
surface $\hat{P}$ in $M(\alpha)$. $P$ is called a presphere if
$\hat{P}$ is a reducing sphere of $M(\alpha)$. Since $M(\alpha)$ is
reducible, the prespheres must exist. Assume $P$ is a presphere such
that $|\partial P|$ is minimal. Then $P$ must be incompressible.

Now suppose $P$ is $\partial$-compressible, with $D$ a
$\partial$-compressing disk. Let $\partial D = u\cup v$, where $u$
is an arc in $P$, and $v$ is an essential arc in $\partial M$. Since
$P$ is incompressible, $v$ is essential on $\partial M-\partial P$.

$\partial$-compressing $P$ along $D$, we get a new surface, which
has one or two new boundary components, depending on whether the two
endpoints of $v$ lie on the different components of $P$. If a new
boundary component is trivial in $\partial M$, we cap off the
component by a disk. In this way, we get a new surface denoted by
$P^\prime$. There are two possibilities:

(1) \ $v$ has endpoints on the different components of $\partial P$.

Now $\hat{P^\prime}$ is also a reducing 2-sphere and $|\partial
P^\prime| < |\partial P|$. It contradicts the assumption that
$|\partial P|$ is minimal.

(2) $v$ has endpoints on the same component of $\partial P$.

$\hat{P^\prime}$ has two components, each of which  is a compressing
disk of $M(\alpha)$, a contradiction. \hfill$\Box$\vskip 5mm

Suppose that $M$ is a simple 3-manifold, and $F$ is a genus at least
two component of $\partial M$. Assume $\alpha$ and $\beta$ are
separating, reducing slopes on $F$. If one of $M(\alpha)$ and
$M(\beta)$, say $M(\beta)$, is $\partial$-reducible, then, by Lemma
4.2 of [SW], $\Delta(\alpha,\beta)=0$. Hence we may assume that
$M(\alpha)$ and $M(\beta)$ are $\partial$-irreducible.

Suppose $\hat{P}$(resp. $\hat{Q}$) is a reducing 2-sphere in
$M(\alpha)$(resp. $M(\beta)$) such that $p=|\partial P|$(resp.
$q=|\partial Q|$) is minimal among all the reducing 2-spheres, where
$P=\hat{P}\cap M$(resp. $Q=\hat{Q}\cap M$). By the proof of Lemma
2.1, $P$ and $Q$ are incompressible and $\partial$-incompressible in
$M$. Isotopy $P$ and $Q$ so that $|P\cap Q|$ is minimal. Then each
component of $P\cap Q$ is either an essential arc or an essential
circle on both $P$ and $Q$.

Let $\Gamma_P$ is a graph in $\hat{P}$ obtained by taking the arc
components of $P \cap Q$ as edges and taking the boundary components
of $P$ as fat vertices. Similarly, we can define $\Gamma_Q$ in
$\hat{Q}$.

{\bf Lemma 2.2.} \ There are no 1-sided disk faces on
$\Gamma_P$(resp.$\Gamma_Q$). \label{lemmano1side} \hfill$\Box$\vskip
5mm

\begin{center}
\includegraphics[totalheight=5.5cm]{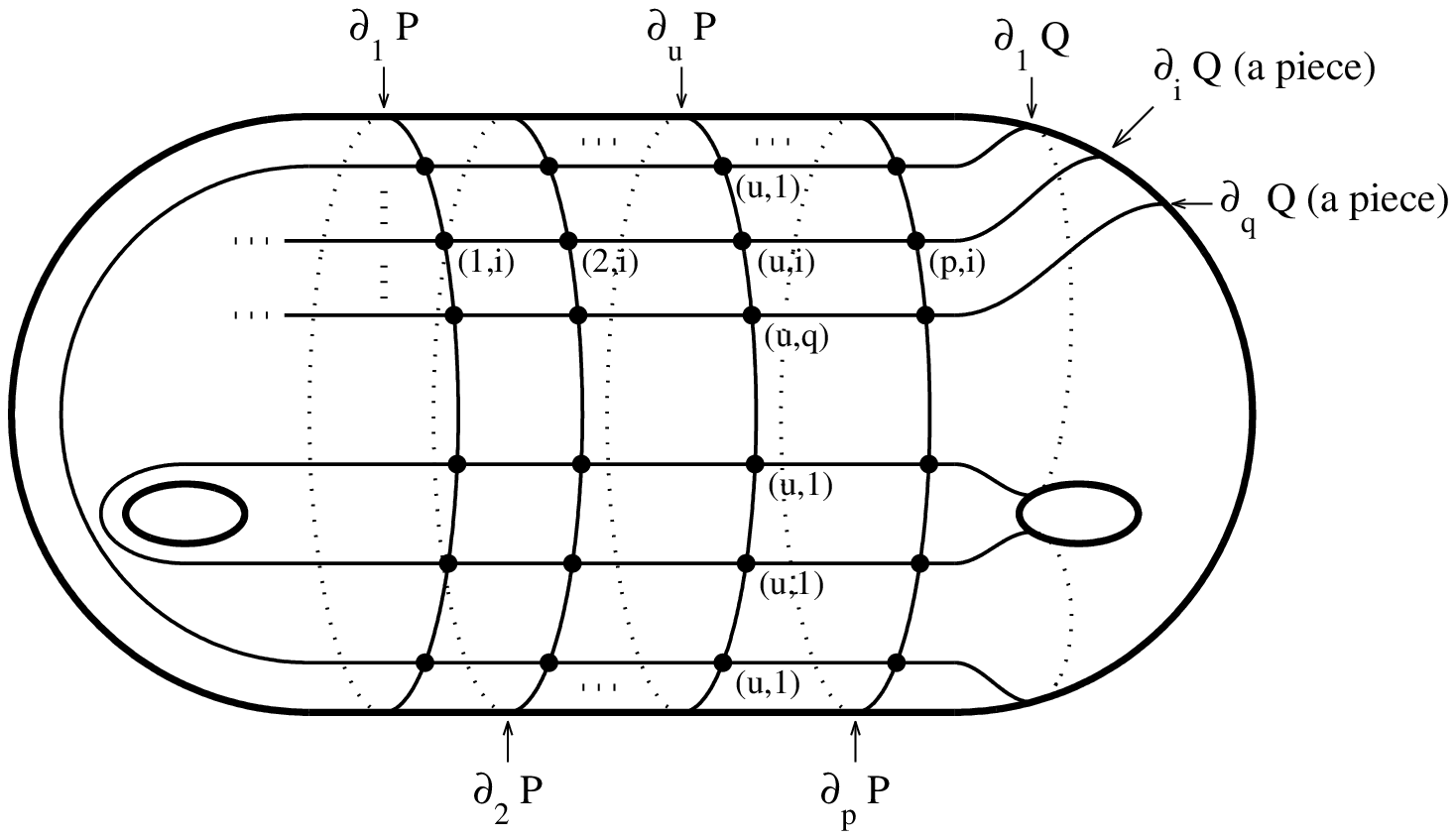}
\begin{center}
Figure 1
\end{center}
\end{center}

\begin{center}
\includegraphics[totalheight=2.5cm]{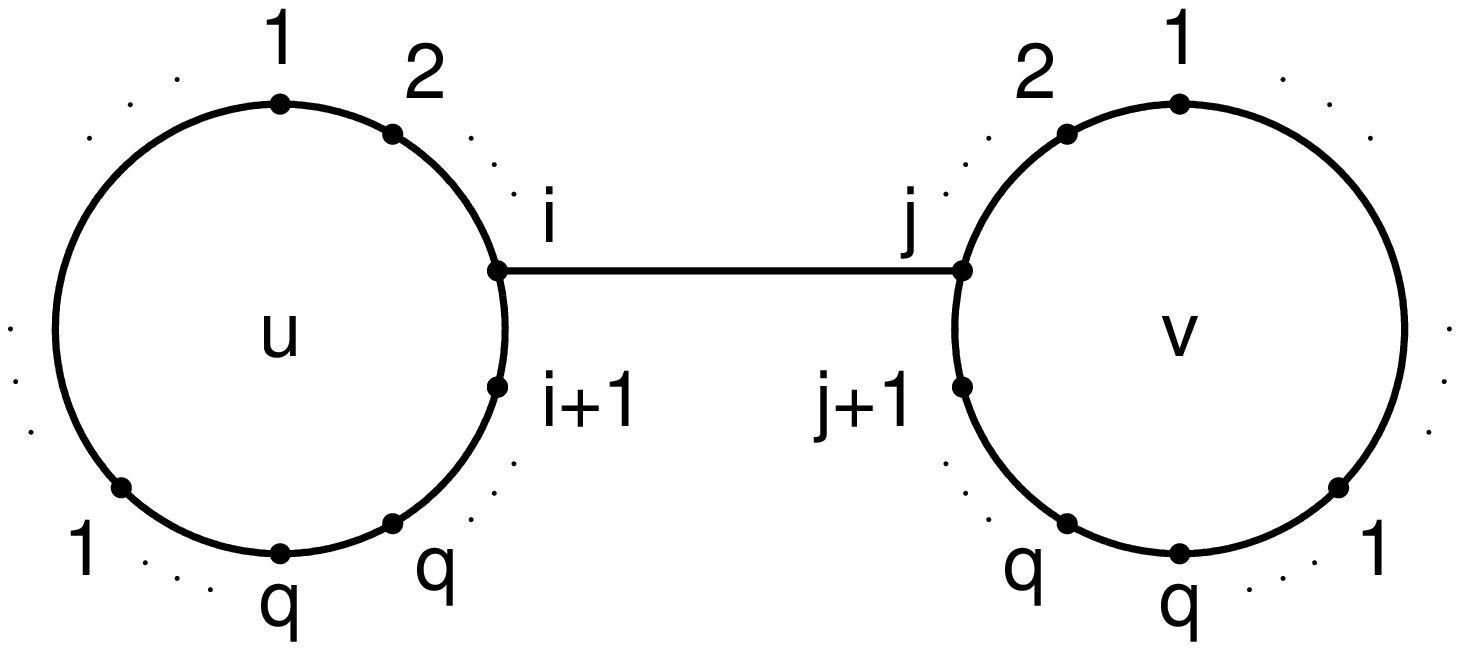}
\includegraphics[totalheight=2.5cm]{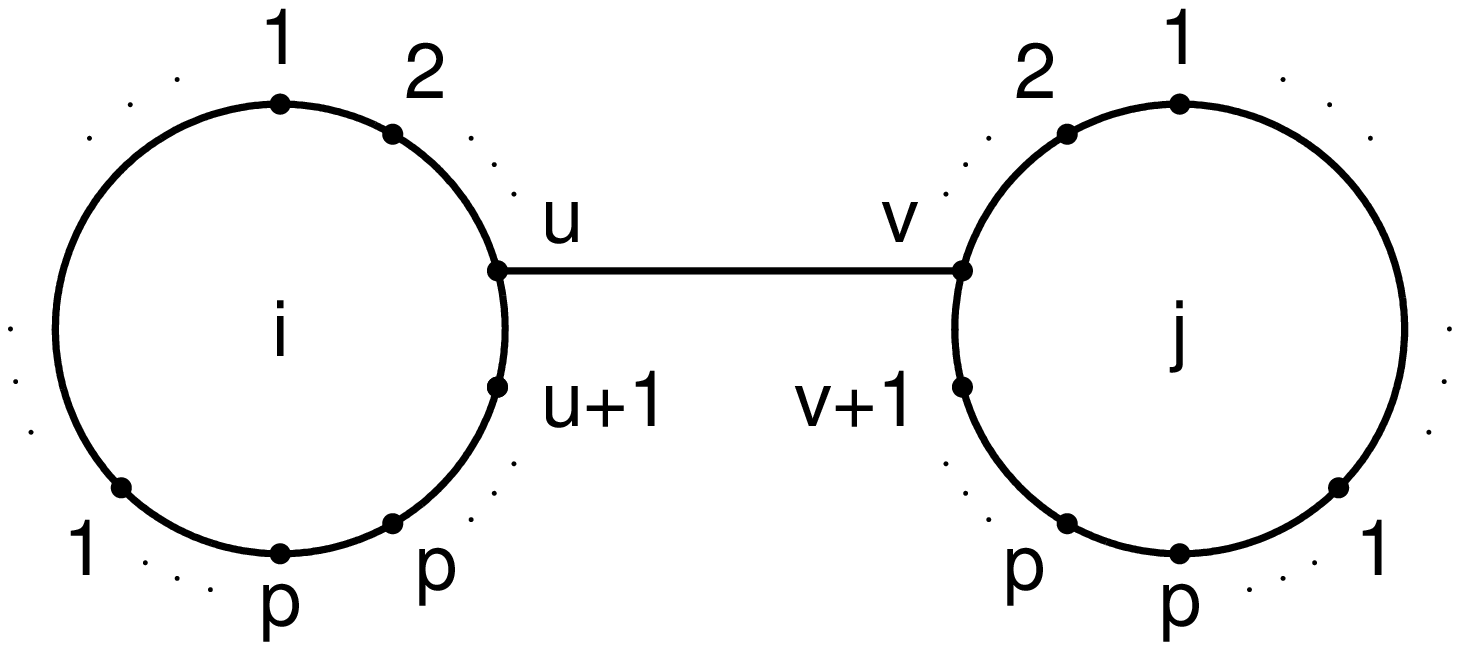}
\begin{center}
Figure 2: Labels on $\Gamma_P$ and $\Gamma_Q$
\end{center}
\end{center}

Number the components of $\partial P$ with $\partial_1 P,
\partial _2 P,\cdots,\partial_u P,\cdots,\partial_{p} P$ consecutively on $\partial
M$, this means that $\partial_u P$ and $\partial_{u+1} P$ bound an
annulus in $\partial M$ with interior disjoint from $P$. See Figure
1.  Similarly, number the components of $\partial Q$ with
$\partial_1 Q,\partial_2 Q,\cdots,\partial_i Q,\cdots,\partial_{q}
Q$. These give  corresponding labels of the vertices of $\Gamma_P$
and $\Gamma_Q$. For an endpoint $x$ of an edge $e$ in $\Gamma_P$, if
it belongs to $\partial_u P\cap
\partial_i Q$, then we label it as $(u,i)$, or $i$(resp. $u$) in $\Gamma_P$(resp.
$\Gamma_Q$) for shortness when $u$(resp. $i$) is specified. See
Figure 2.  Now each edge $e$ of $\Gamma_P$ has been labeled with
$(u,i)-(v,j)$, or $i-j$(resp. $u-v$) in $\Gamma_P$(resp. $\Gamma_Q$)
for shortness. See Figure 2.  When we travel around $\partial_u P$,
the labels appear in the order $1, 2,\cdots,q, q, \cdots, 2,
1,\cdots$(repeated $\Delta(\alpha,\beta)/2$ times). Note that
$\Gamma_Q$ have the same property.

\section {Parity rule}
\label{parityrule} We first sign the endpoints of the edges in
$\Gamma_P$(and in $\Gamma_Q$). Fix the directions on $\alpha$ and
$\beta$. Then each point in $\alpha \cap \beta$ can be signed
``$+$'' or ``$-$'' depending on whether the direction determined by
right-hand rule from $\alpha$ to $\beta$ is pointed to the outside
of $M$ or to the inside of $M$. See Figure 3. Since $\alpha$ and
$\beta$ is separating, the signs ``$+$'' and ``$-$'' appear
alternately on both $\alpha$ and $\beta$.
\begin{center}
\includegraphics[totalheight=8cm]{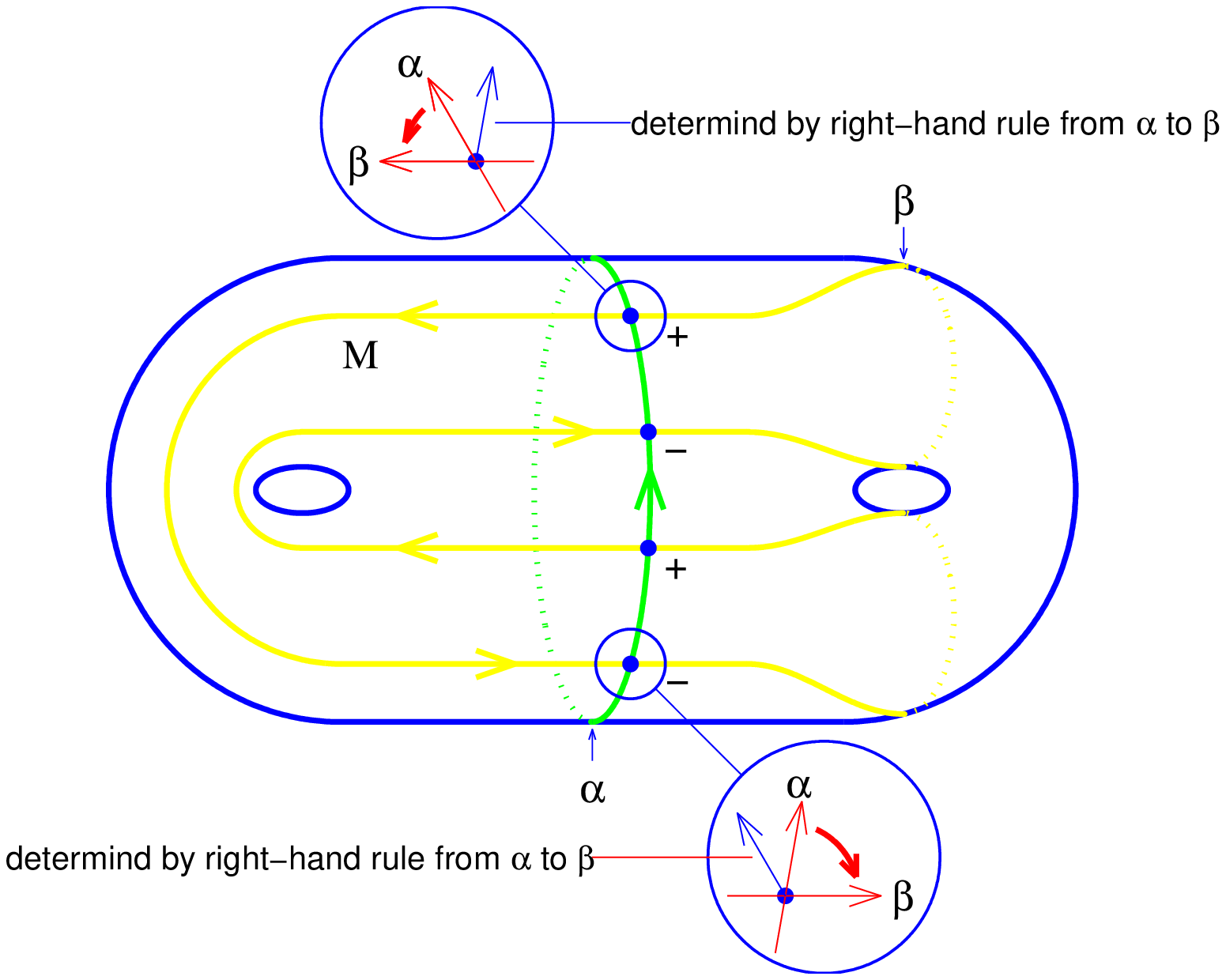}
\begin{center}
Figure 3
\end{center}
\end{center}

\begin{center}
\includegraphics[totalheight=6cm]{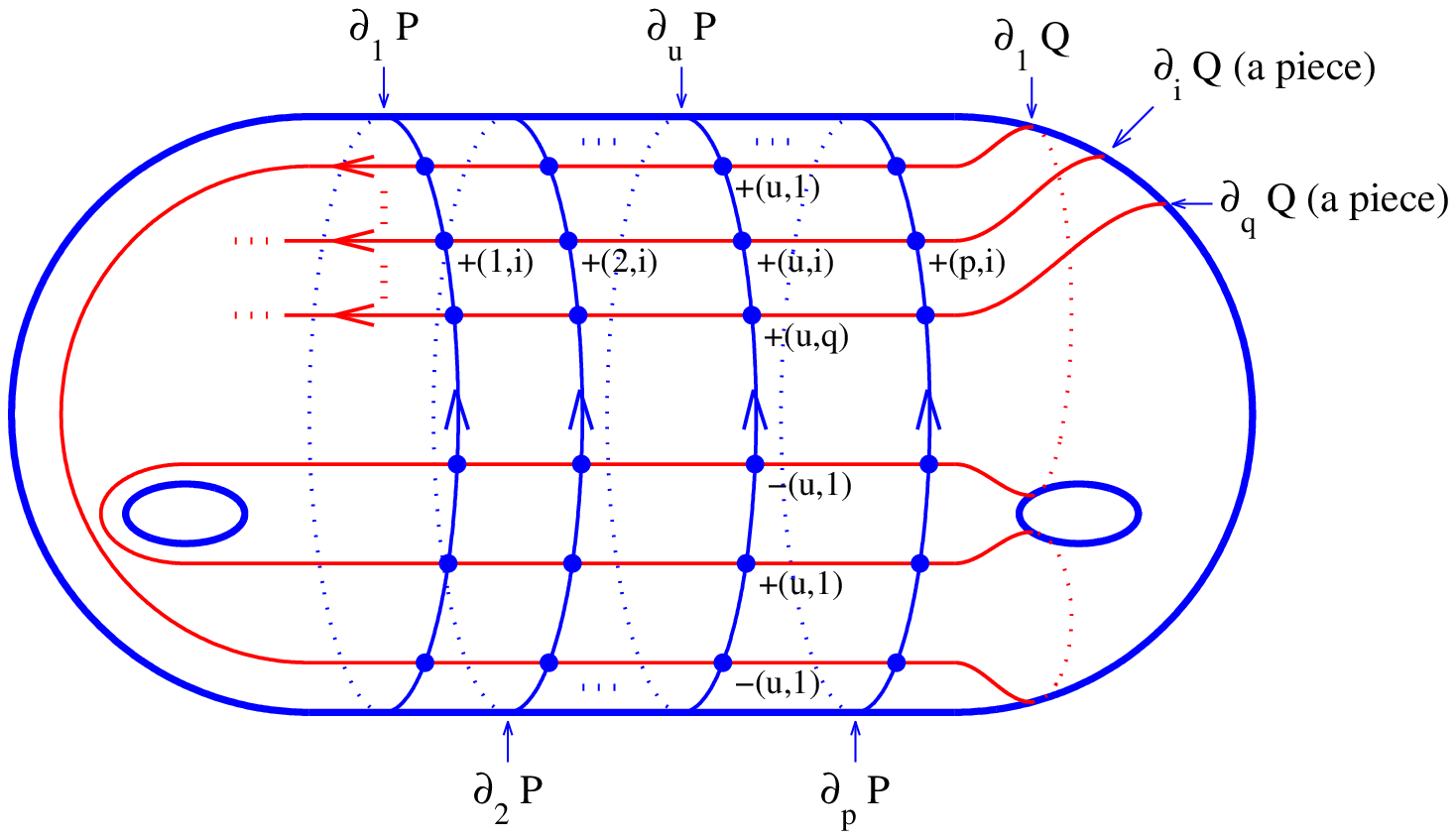}
\begin{center}
Figure 4: Signs on $\partial P\cap\partial Q$
\end{center}
\end{center}

Give a direction to each boundary components $\partial P$(resp.
$\partial Q$) such that they are all parallel to $\alpha$(resp.
$\beta$) on $\partial M$. Then each point $x\in \partial
P\cap\partial Q$ can be signed as above. We denoted by $c(x)$ the
sign of $x$. See Figure 4. Now the signed labels appear on
$\partial_u P$ as $+1,+2,\cdots,+q,-q,\cdots,-2,-1,\cdots$,
(repeated $\Delta(\alpha,\beta)/2$ times). See Figure 5.

\begin{center}
\includegraphics[totalheight=2.5cm]{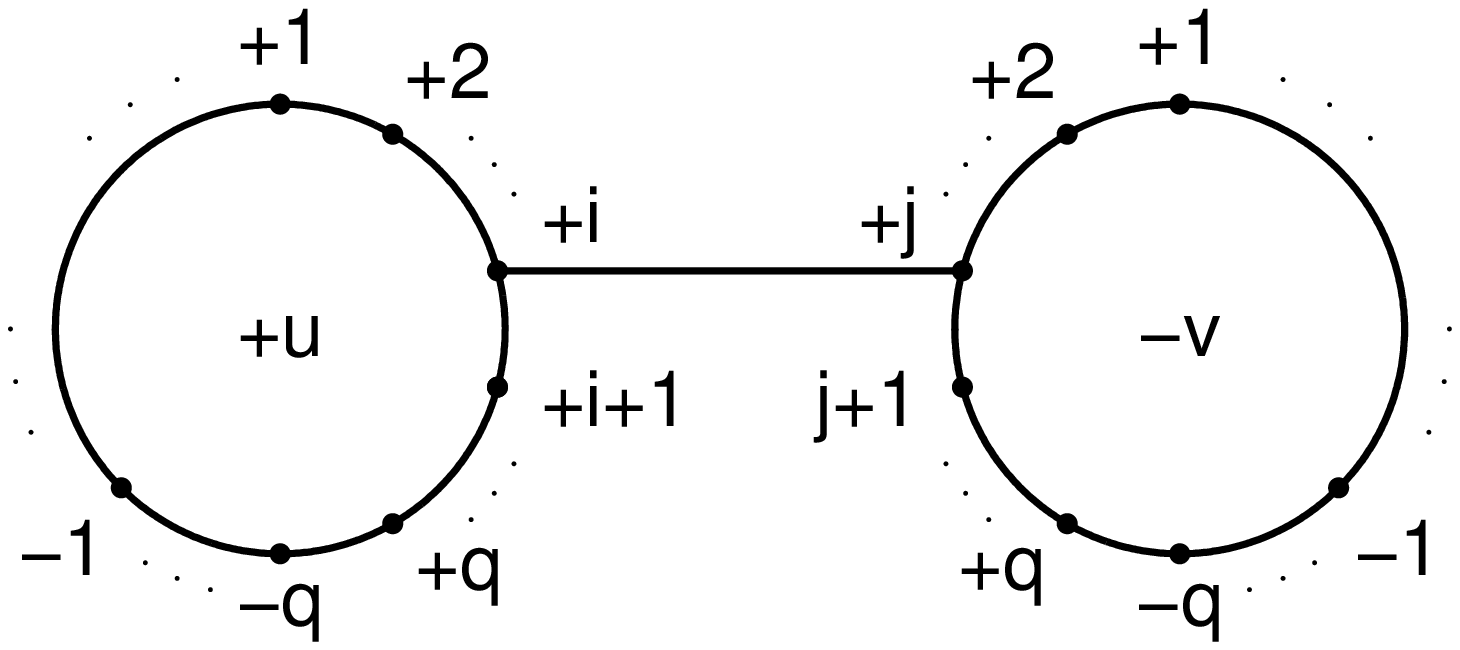}
\includegraphics[totalheight=2.5cm]{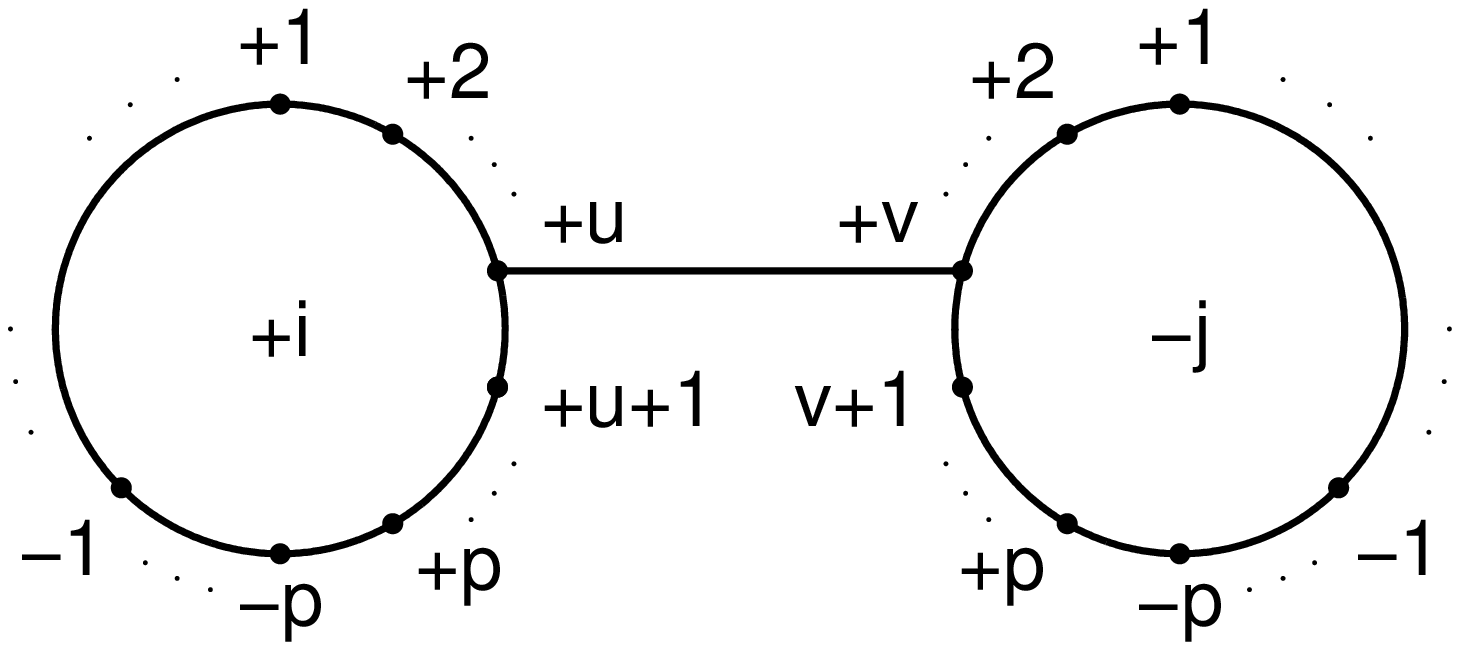}
\begin{center}
Figure 5
\end{center}
\end{center}

Now we sign the vertices of $\Gamma_P$. Suppose $P\times [0,1]$ be a
thin regular neighborhood of $P$ in $M$. Let $P^+=P\times 1$ and
$P^-=P\times 0$. For some $1\leq u\leq p, 1\leq i\leq q$, let $c$ be
a component of $\partial_u P\times[0,1] \cap \partial_i Q$ with the
induced direction of $\partial_i Q$. We define the sign of
$\partial_u P$ as follows:

(1) Suppose $c$ intersects $\partial_u P$ at a ``$+$'' point, we
define the sign of $\partial_u P$ is ``$+$''(resp. ``$-$'') if the
direction of $c$ is from $P^+$ to $P^-$(resp. from $P^-$ to $P^+$).

(2) Suppose $c$ intersects $\partial_u P$ at a ``$-$'' point, we
define the sign of $\partial_u P$ is ``$+$''(resp. ``$-$'') if the
direction of $c$ is from $P^-$ to $P^+$(resp. from $P^+$ to $P^-$).

Since each component of $\partial Q$ has the same direction with
$\beta$ on $F$, the definition as above is independent of the
choices of $c$ and $i$.

For example, in Figure 6 and Figure 7, the signs of $\partial_u P$,
$\partial_v P$ and $\partial_w P$ are ``$+$'', ``$-$'' and ``$-$''
respectively.

Since $M$ is orientable,  $\partial_u P$ and $\partial_v P$ have the
same direction on $P$ when $\partial_u P$ and $\partial_v P$ have
the same signs. This means the labels $+1,+2,\cdots,+q$,
$-q,\cdots,-1$ of the edge-endpoints appear on both $\partial_u P$
and $\partial_v P$ are in the same direction in $\Gamma_P$.
Similarly, the labels $+1,+2,\cdots,+q,-q\cdots,-1$ appear in
opposite the directions when $\partial_u P$ and $\partial_v P$ have
different signs. See Figure 7.

We may define the sign $s(i)$ of $\partial_i Q$ in $\Gamma_Q$.

\begin{center}
\includegraphics[totalheight=6cm]{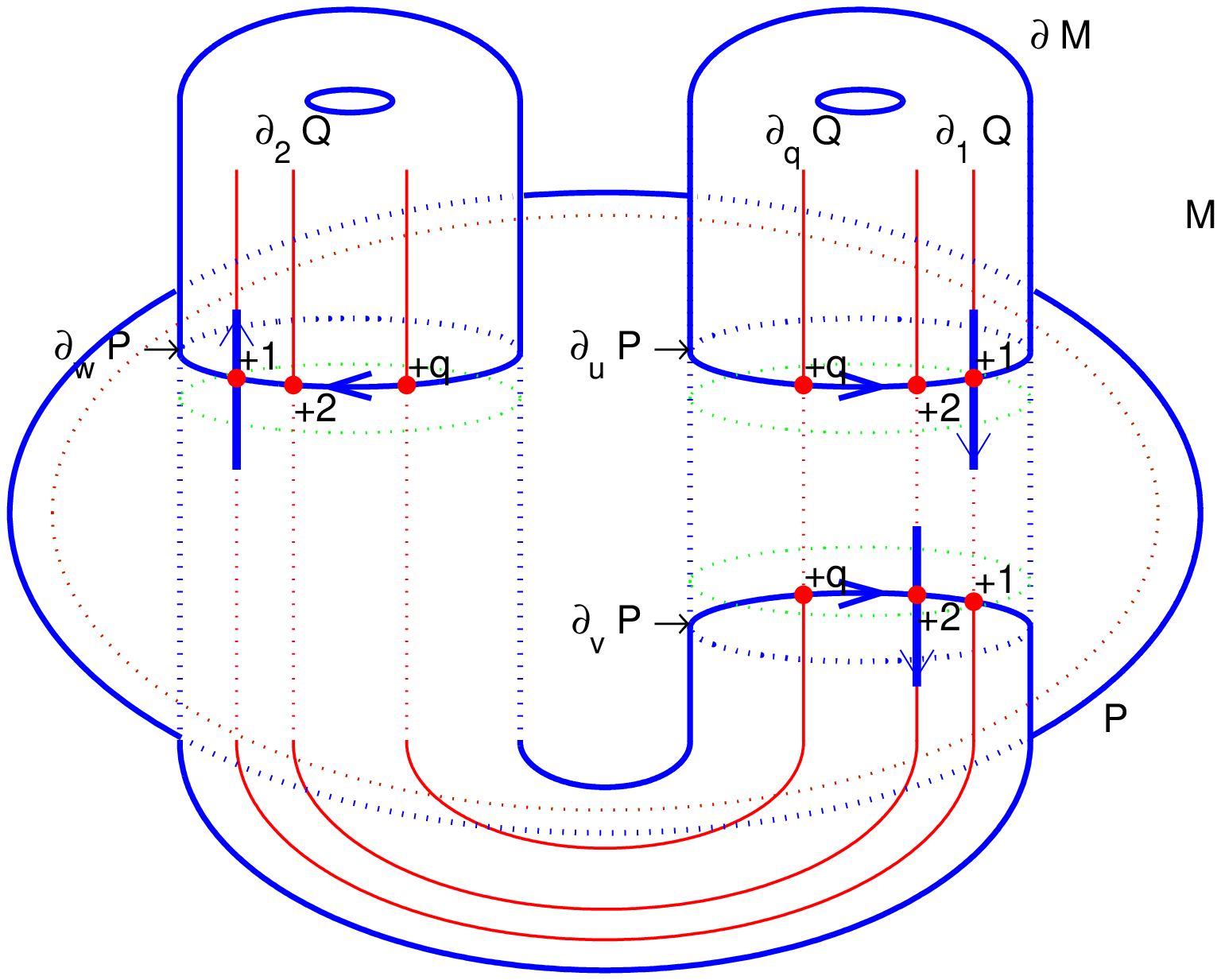}
\begin{center}
Figure 6: Signs on $\partial P$
\end{center}
\end{center}

\begin{center}
\includegraphics[totalheight=4cm]{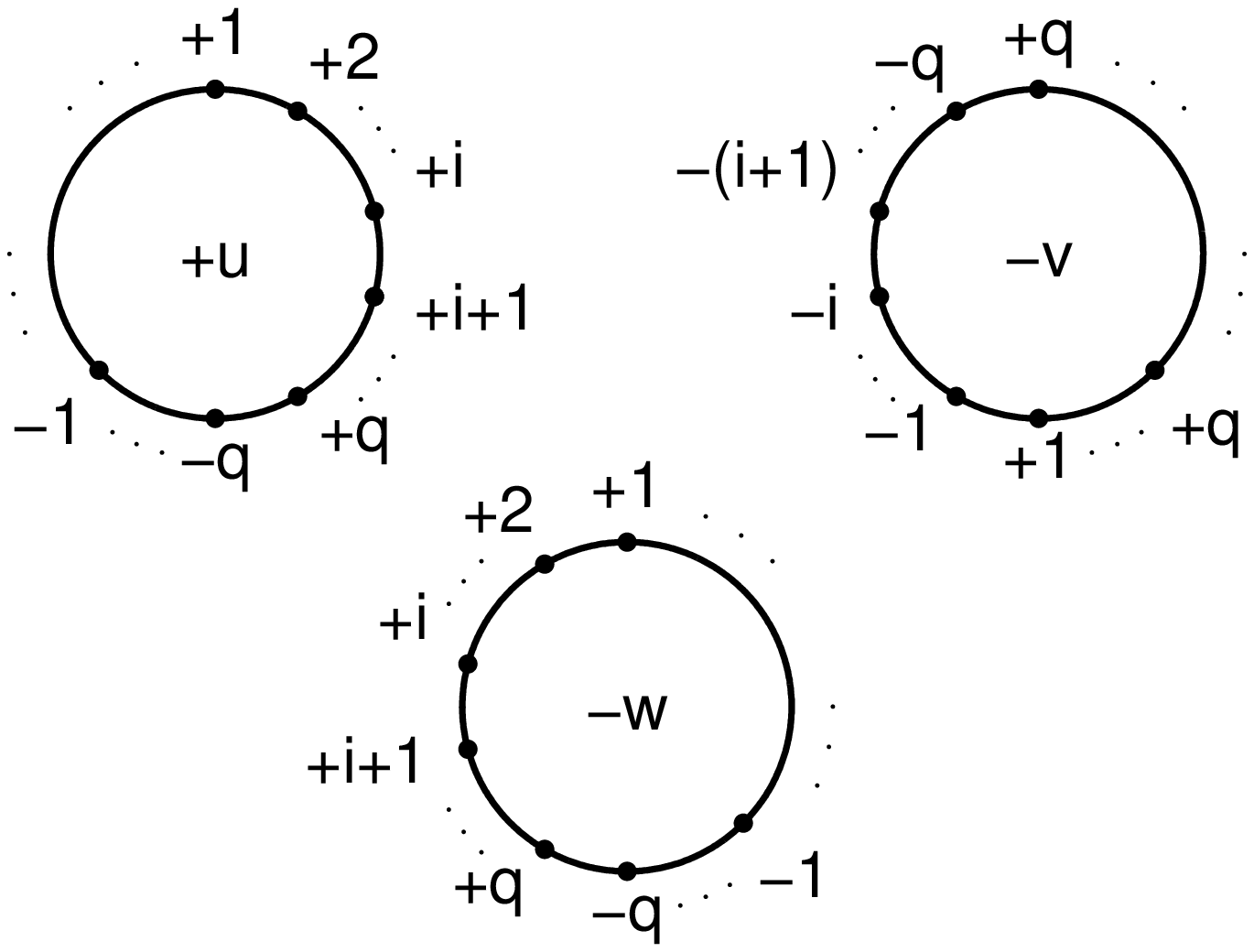}
\begin{center}
Figure 7: Signs on $\partial P$
\end{center}
\end{center}

The labels with the signs defined as above are said to be Type A.
Now we have Parity rule A:

{\bf Lemma 3.1(Parity rule A).} \ For an edge $e$ in $\Gamma_P$(and
$\Gamma_Q$) with its endpoints $x$ labeled $(u,i)$ and $y$ labeled
$(v,j)$, the following equality holds:

\begin{center}
$s(i)s(j)s(u)s(v)c(x)c(y)=-1$ ~~~~~ ($*$).
\end{center}\label{lemmaparityruleA}

{\bf Proof.} \ Let $P\times I$ be a thin regular neighborhood of $P$
in $M$. Then $e\times I\subset Q$ and $x\times I\subset \partial_u
P\times I$ and $y\times I\subset
\partial_v P\times I$.

Now there are four possibilities:

{\bf Case 1.} \ $s(i)=s(j)$ and $c(x)=c(y)$ as in Figure 8(a).

Since $s(i)=s(j)$, $\partial_i Q$ and $\partial_j Q$ have the same
direction. In this case, $x\times I$ and $y\times I$ have the
opposite directions(as in Figure 8(a). Since $c(x)=c(y)$, by the
definitions of $s(u)$ and $s(v)$, $s(u)\neq s(v)$. Hence the
equality ($*$) holds.

{\bf Case 2} \ $s(i)=s(j)$ and $c(x)\neq c(y)$ as in Figure 8(b).

Since $s(i)=s(j)$, $\partial_i Q$ and $\partial_j Q$ have the same
direction. In this case, $x\times I$ and $y\times I$ have the
opposite directions as in Figure 8(b). Since $c(x)\neq c(y)$, by the
definitions of $s(u)$ and $s(v)$, $s(u)=s(v)$. Hence the equality
($*$) holds.

{\bf Case 3} \ $s(i)\neq s(j)$ and $c(x)= c(y)$ as in Figure 8(c). .

Since $s(i)\neq s(j)$, $\partial_i Q$ and $\partial_j Q$ have
opposite directions. In this case, $x\times I$ and $y\times I$ have
the same direction as in Figure 8(c). Since $c(x)=c(y)$, by the
definitions of $s(u)$ and $s(v)$, $s(u)=s(v)$. Hence the equality
($*$) holds.

{\bf Case 4} \ $s(i)\neq s(j)$ and $c(x)\neq c(y)$ as in Figure
8(d).

Since $s(i)\neq s(j)$, $\partial_i Q$ and $\partial_j Q$ have
opposite directions. In this case, $x\times I$ and $y\times I$ have
the same direction as in Figure 8(d).  Since $c(x)\neq(y)$, by the
definitions of $s(u)$ and $s(v)$, $s(u)\neq s(v)$.

Hence the equality ($*$) holds. \hfill$\Box$\vskip 5mm

\begin{center}

\includegraphics[totalheight=2.8cm]{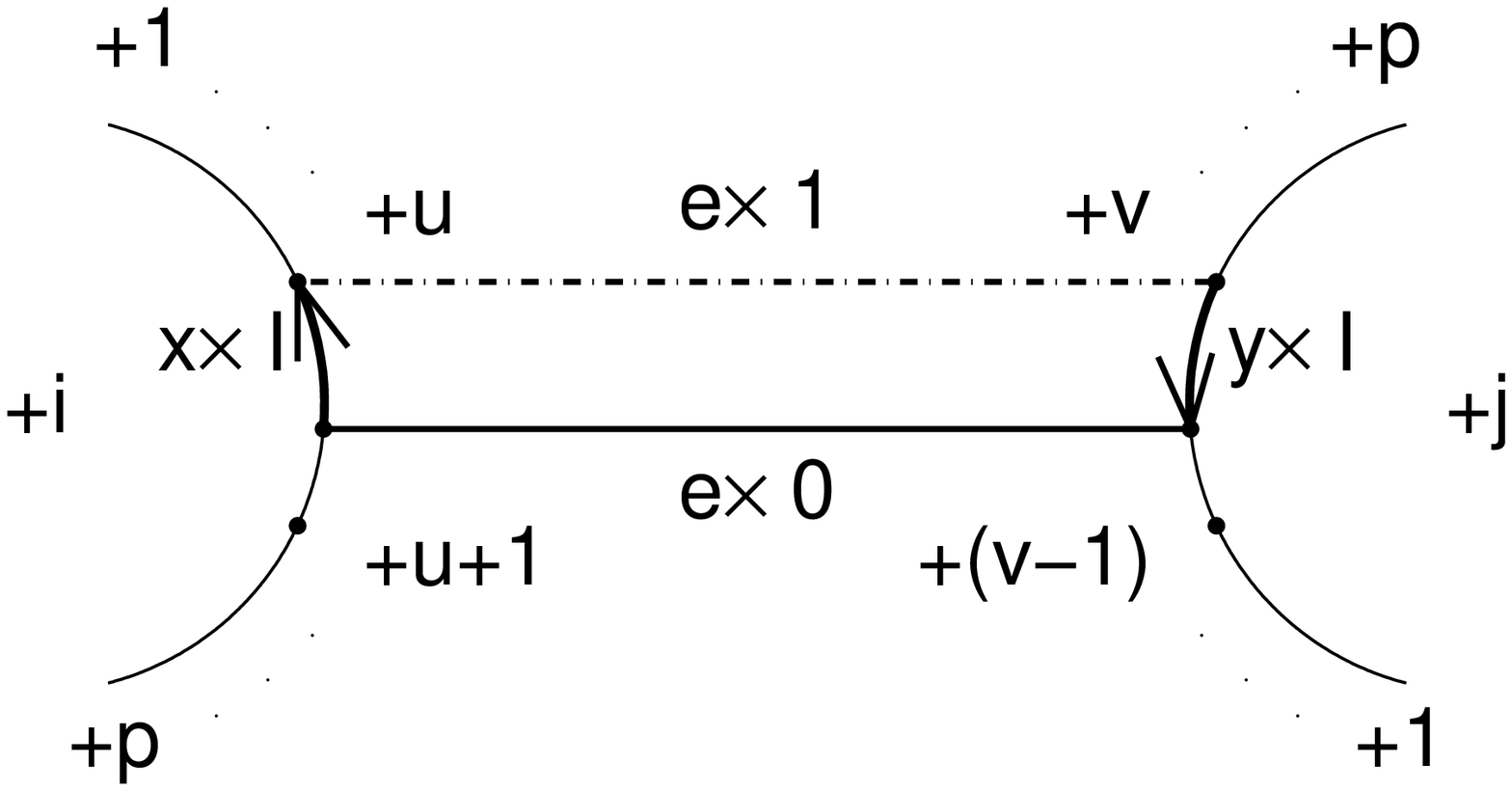}\hfill
\includegraphics[totalheight=2.8cm]{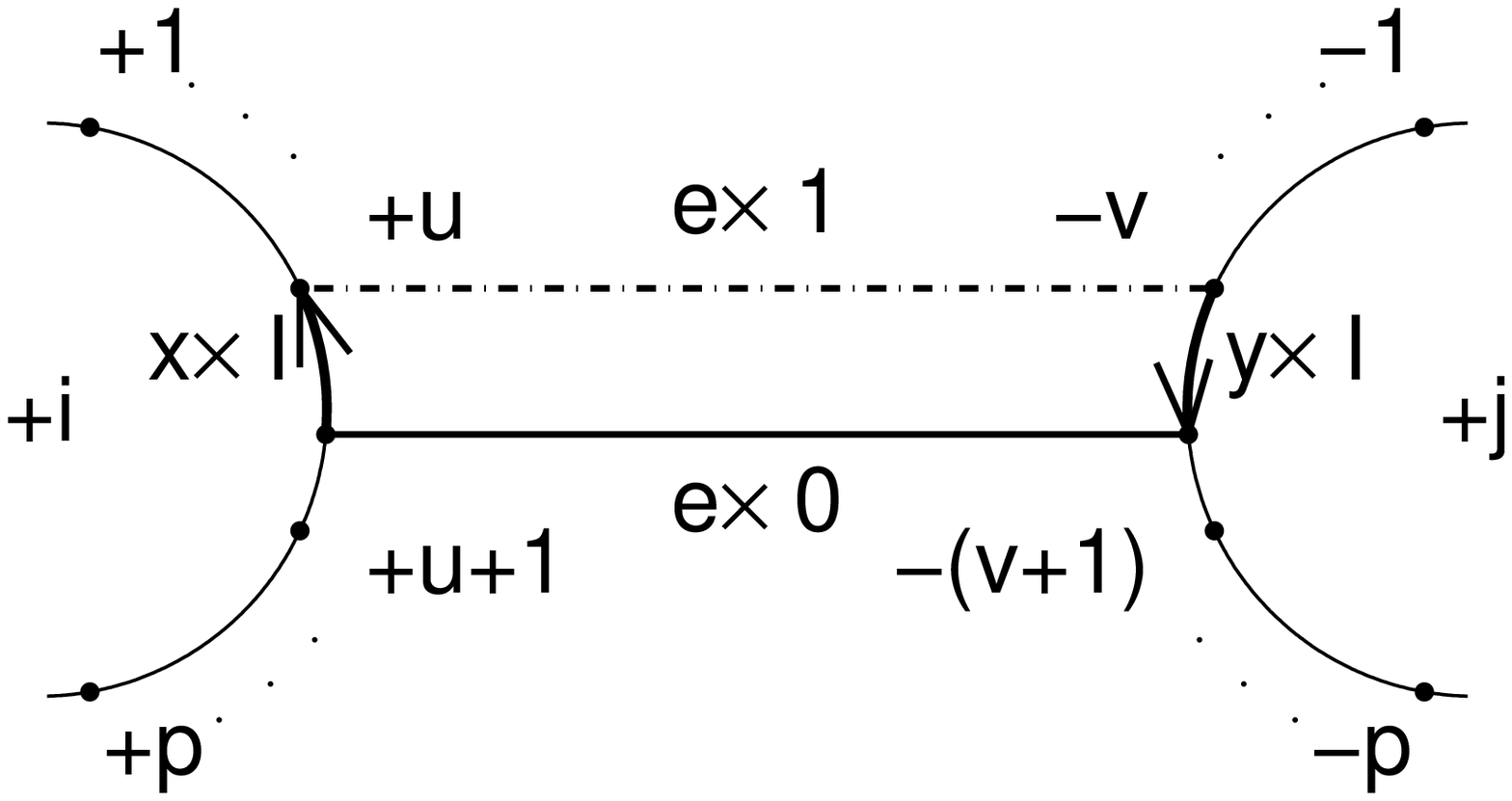}
\begin{center}(a)~~~~~~~~~~~~~~~~~~~~~~~~~~~~~~~~~~~~~~~~~~~~~~~~~~~~~~~~~~~~~~~~~(b)
\end{center}

\includegraphics[totalheight=2.8cm]{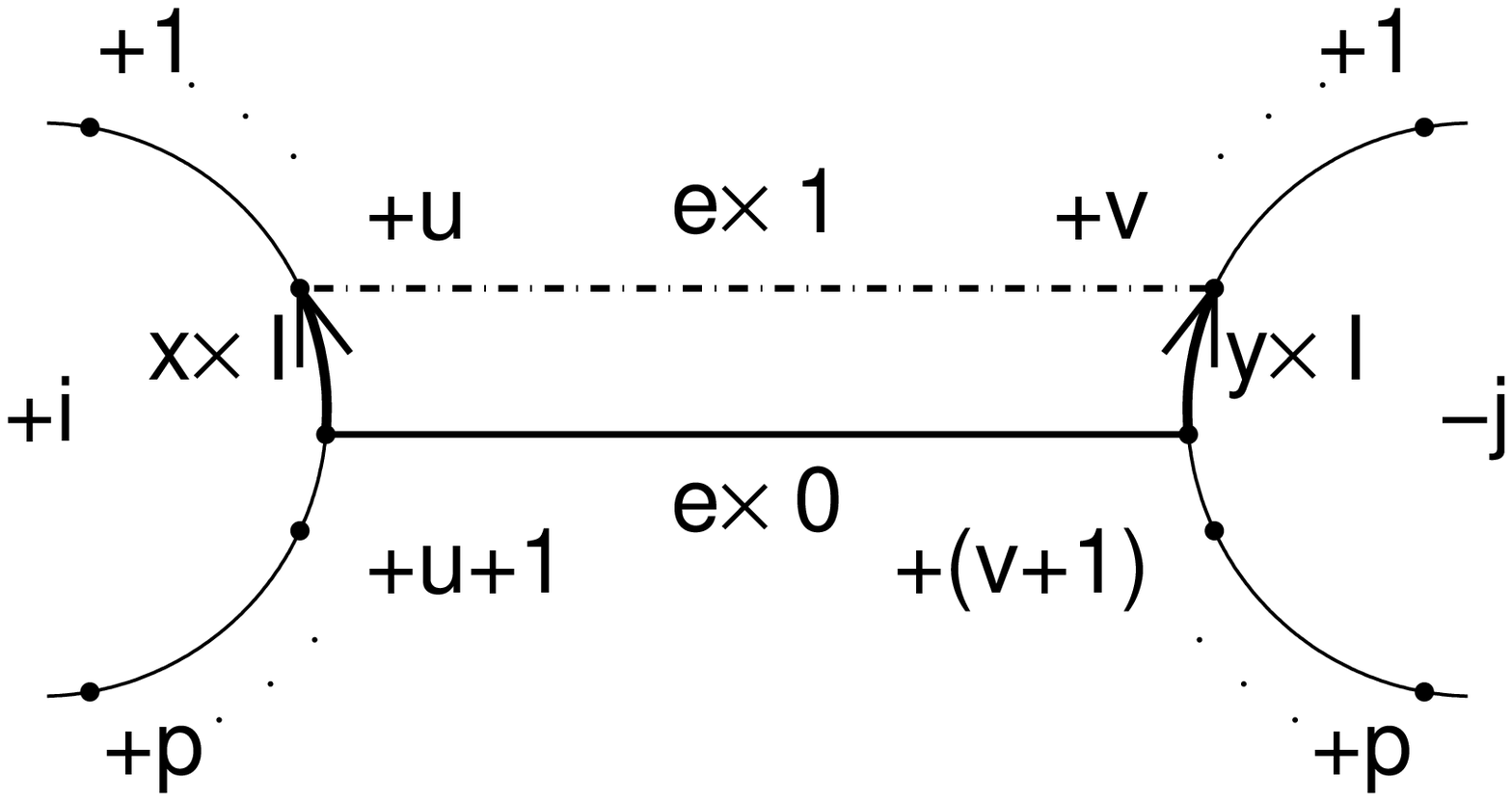}\hfill
\includegraphics[totalheight=2.8cm]{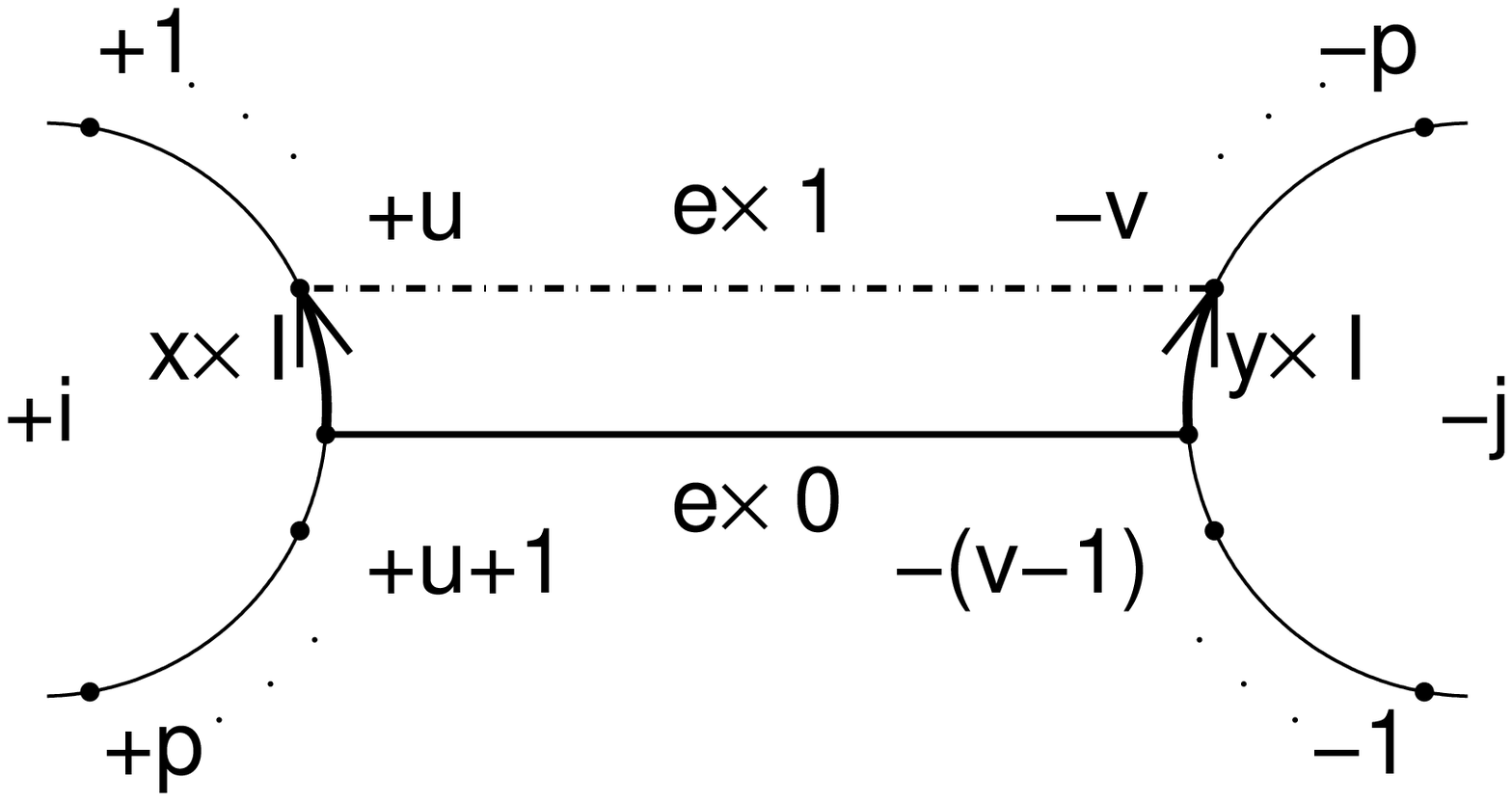}
\begin{center}(c)~~~~~~~~~~~~~~~~~~~~~~~~~~~~~~~~~~~~~~~~~~~~~~~~~~~~~~~~~~~~~~~~~(d)
\end{center}
\begin{center}
Figure 8
\end{center}

\end{center}

Now suppose that $e$ is an edge of $\Gamma_P$ with $\partial e=x\cup
y$, and $x$ is labeled $(u,i)$. Let $g(x)=c(x)\times s(u)$. Then the
signed label $g(x)i$ of $x$ is said to be Type B.

{\bf Remark ($*$) } \ Under Type B labels, the labels $+1$,
$+2\cdots$, $+q,-q,\cdots$ appear in the same direction on all the
vertices of $\Gamma_P$. See Figure 9.

\begin{center}

\includegraphics[totalheight=4.5cm]{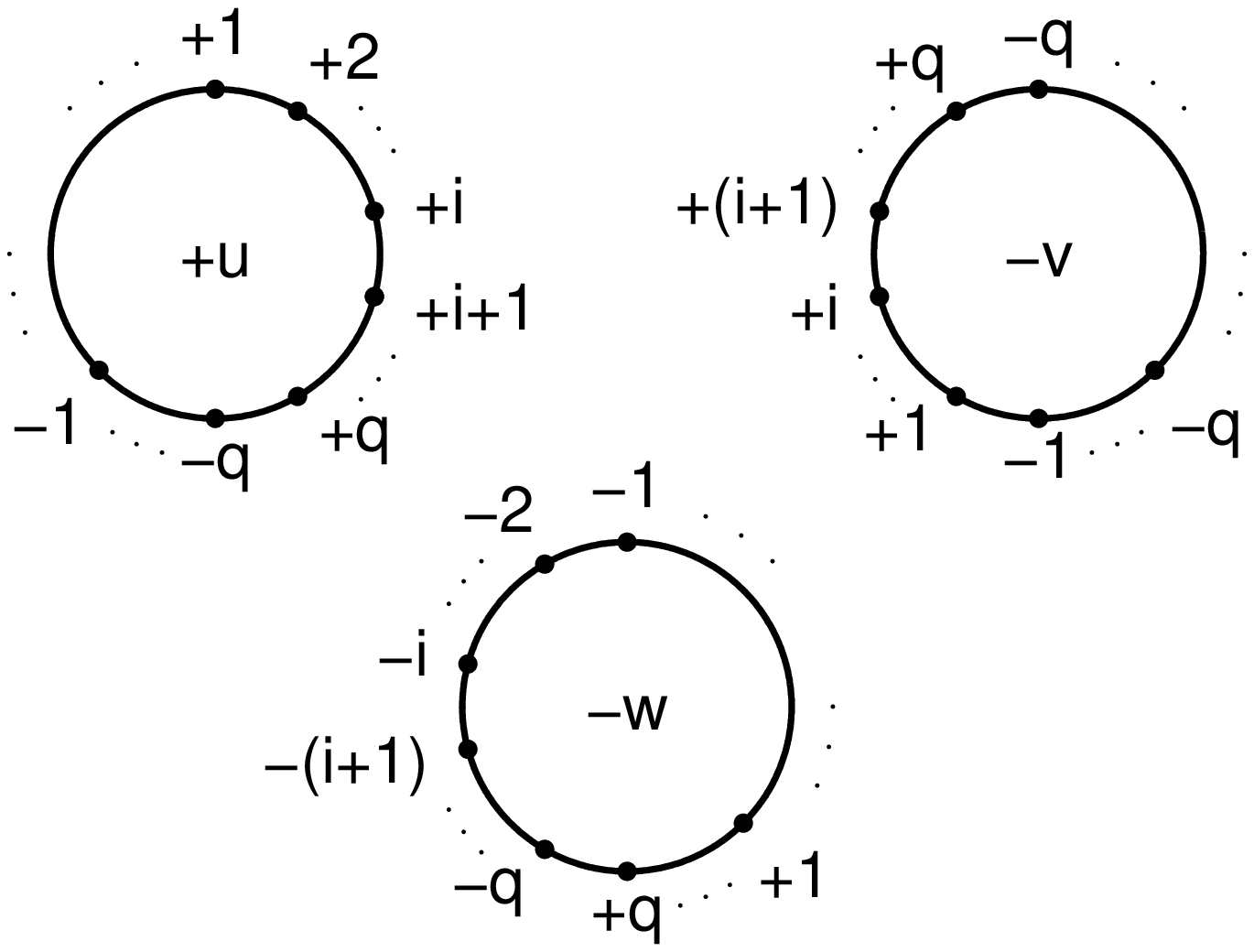}
\begin{center}
Figure 9: Type B labels
\end{center}
\end{center}

By Lemma 3.1, we have the parity rule for Type B labels.

{\bf Lemma 3.2(Parity rule B).} \ Let $e$ be an edge $e$ in
$\Gamma_P$ with its endpoints $x$ labeled $(u,i)$ and $y$ labeled
$(v,j)$, then $s(i)s(j)g(x)g(y)=-1$.\label{lemmaparityruleB}
\hfill$\Box$\vskip 5mm

{\bf Lemma 3.3.} \ Let $e$ be an edge $e$ in $\Gamma_P$ with its
endpoints $x$ labeled $(u,i)$ and $y$ labeled $(v,i)$. Then
$g(x)\neq g(y)$. \label{lemmaparityrule} \hfill$\Box$\vskip 5mm

\section {S-cycles}
\label{scharlemanncycles}

In this section, the definitions of a cycle, the length of a cycle,
a disk face and parallel edges are standard, see [GL1], [SW] and
[W].

\begin{center}
\centering
\includegraphics[totalheight=3.5cm]{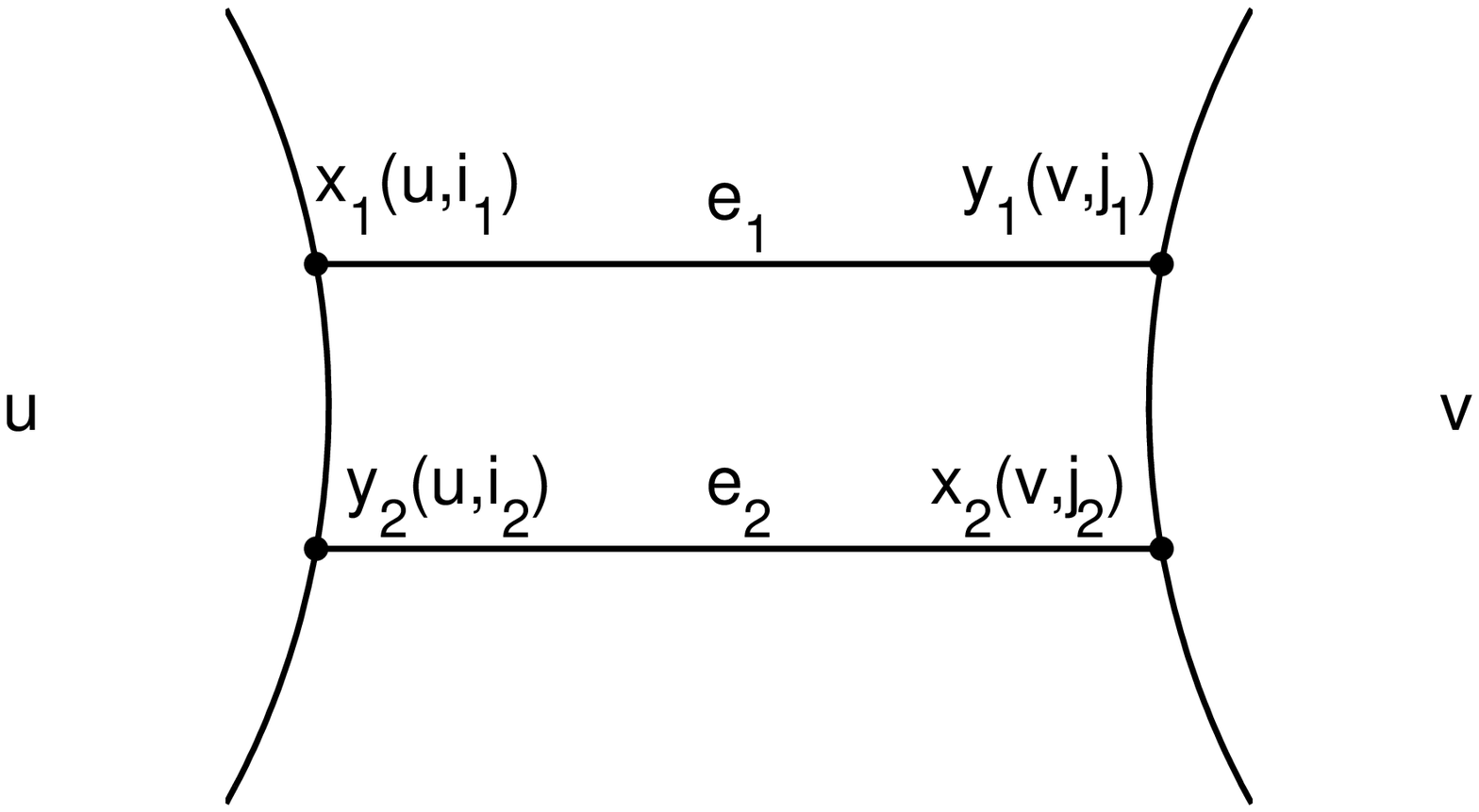}

Figure 10:
\end{center}

Suppose a length two cycle $C=\{e_1,e_2\}$ bounds a disk-face in
$\Gamma_P$, where $\partial e_1=x_1\cup y_1$ with $x_1$ labeled
$(u,i_1)$ and $y_1$ labeled $(v,j_1)$, and $\partial e_2=x_{2}\cup
y_2$ with $x_2$ labeled $(v,j_2)$ and $y_2$ labeled $(u,i_2)$. See
Figure 10. $C$ is said to be a virtual S-cycle if
$g(x_1)i_1=g(x_2)j_2$ and $g(y_2)i_2=g(y_1)j_1$. In this case,
$\{i_1,j_1\}$ is called the label pair of $C$. Furthermore if
$i_1\neq j_1$, then $C$ is called an S-cycle.

\begin{center}

\includegraphics[totalheight=3cm]{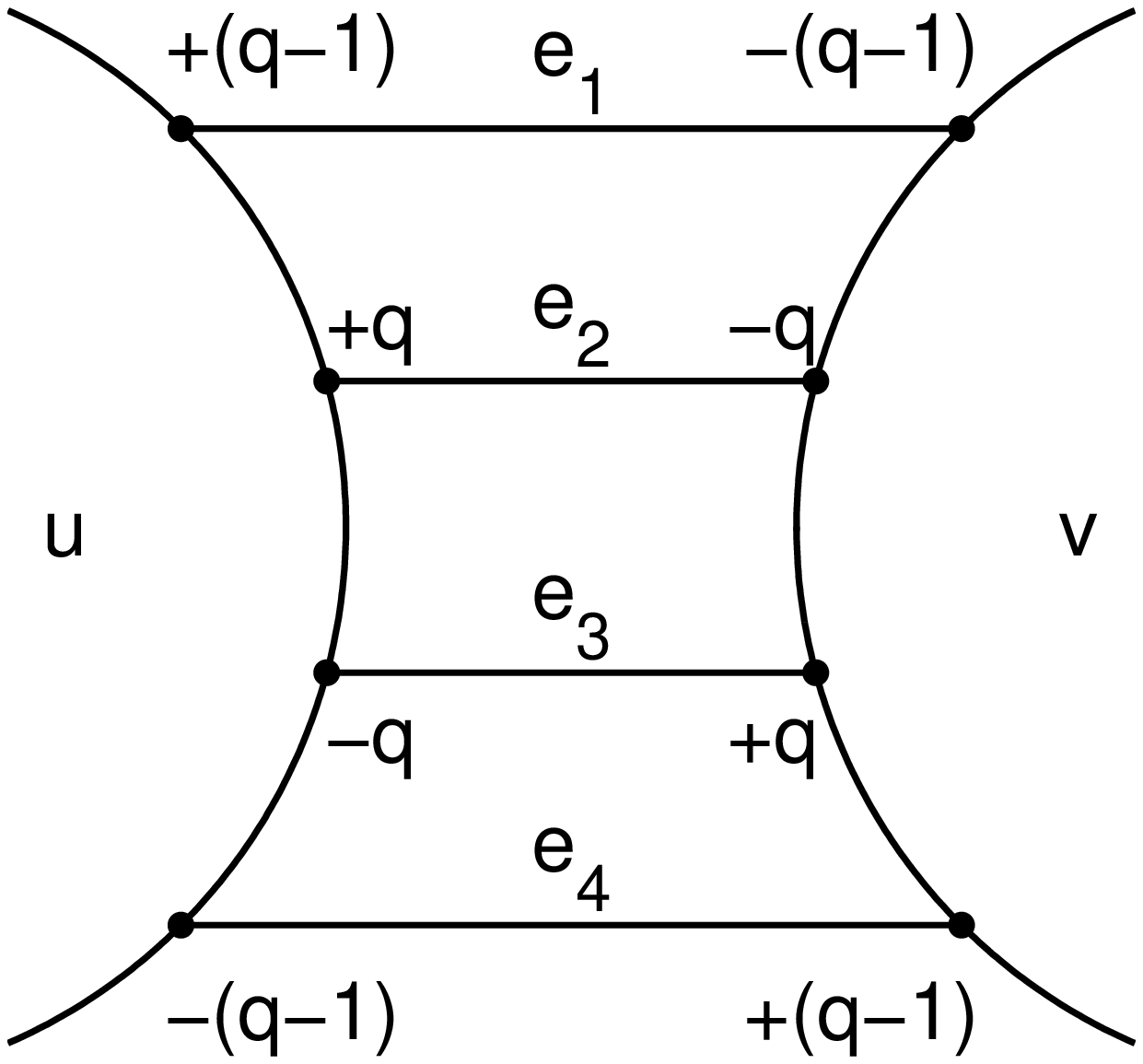}\hfill
\includegraphics[totalheight=3cm]{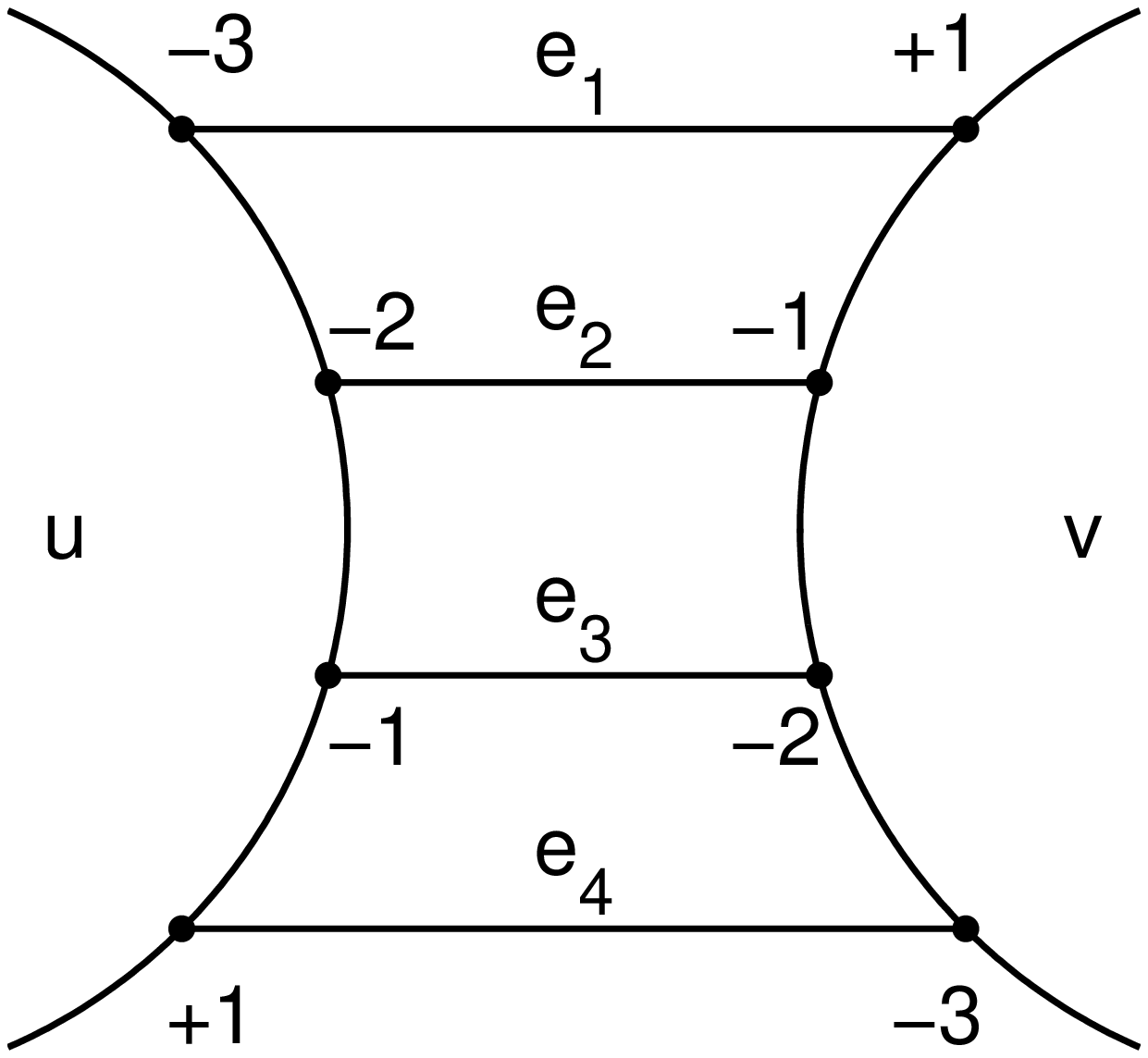}\hfill
\includegraphics[totalheight=3cm]{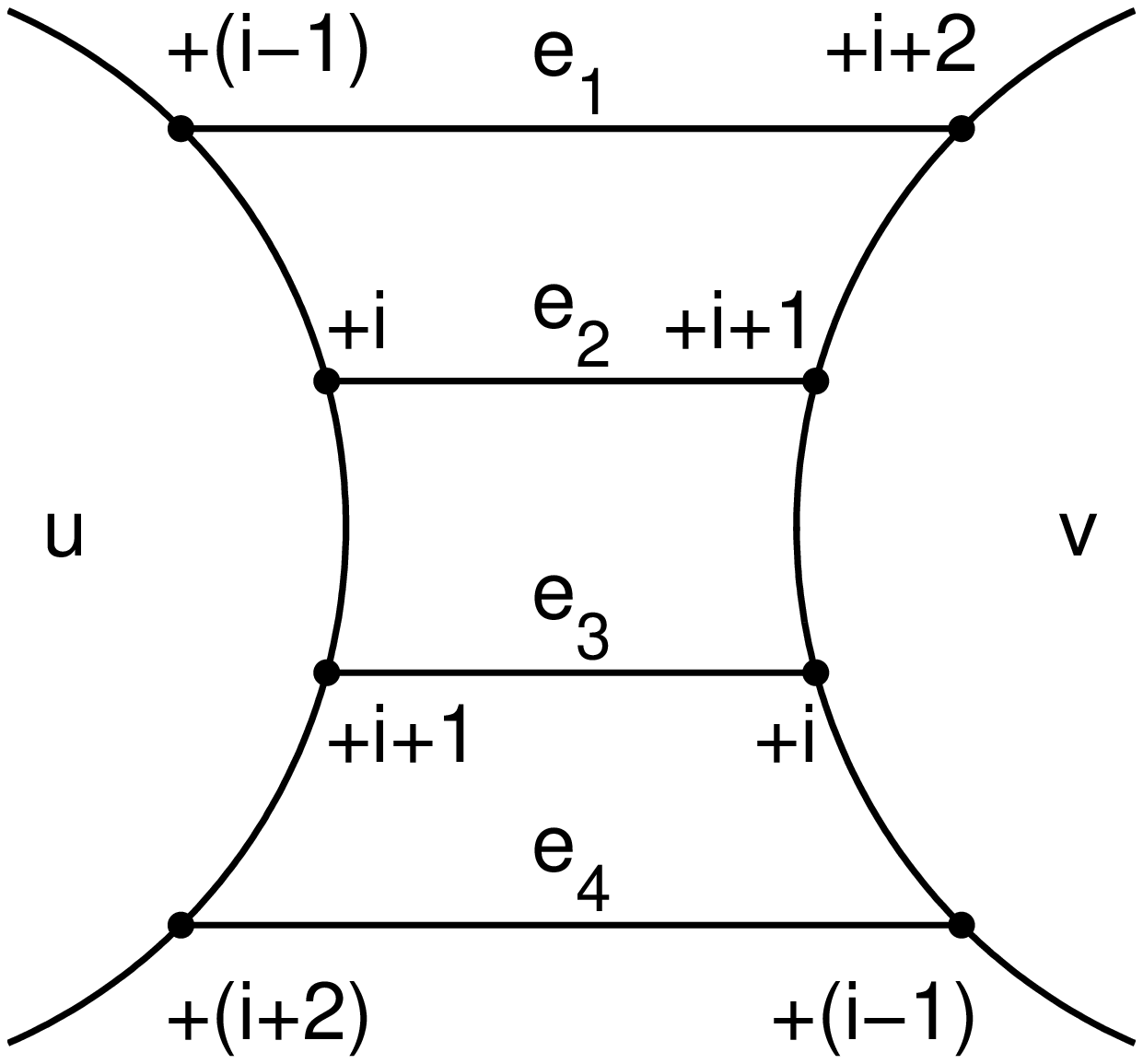}

\begin{center}
(a)~~~~~~~~~~~~~~~~~~~~~~~~~~~~~~~~~~~~~~~~~~~~(b)~~~~~~~~~~~~~~~~~~~~~~~~~~~~~~~~~~~~~~~~~~~(c)
\end{center}
\begin{center}
Figure 11
\end{center}
\end{center}

{\bf Lemma 4.1.} \ A virtual S-cycle is either an S-cycle, or its
label pair is one of $\{1,1\}$ and $\{q, q\}$.

{\bf Proof} \ Let $\{e_1,e_2\}$ be an S-cycle defined as above. If
$i_1\neq i_2$, then it is an S-cycle. If $i_1=i_2$, then
$i_1=j_1=i_2=j_2$. Hence either $i_1=1$ or $i_1=q$.
\hfill$\Box$\vskip 5mm

A set of four adjacent parallel edges, say $\{e_1, e_2, e_3, e_4\}$,
in $\Gamma_P$ is called a virtual extended S-cycle if $\{e_2, e_3\}$
is an S-cycle.

A virtual extended S-cycle $\{e_1, e_2, e_3, e_4\}$ is called an
extended S-cycle if $\{e_2,e_3\}$ is not an S-cycle labeled
$\{1,2\}$ or $\{q,q-1\}$.

For examples, in Figure 11(a), $\{e_2,e_3\}$ is a virtual S-cycle
rather than an S-cycle, and $\{e_1, e_2, e_3, e_4\}$ is a virtual
extended S-cycle rather than an extended S-cycle; in Figure 11(b),
$\{e_2, e_3\}$ is an S-cycle, but $\{e_1, e_2, e_3, e_4\}$ is a
virtual extended S-cycle rather than an extended S-cycle; in Figure
11(c), $\{e_1, e_2, e_3, e_4\}$ is an extended S-cycle.

{\bf Lemma 4.2.} \ (1) $\Gamma_P$ can not contain two S-cycles with
distinct label pairs. \label{lemmadistingctlabelpair}

(2) $\Gamma_P$ contains no extended S-cycles.
\label{lemmaextededcycle}

{\bf Proof} \ The proof follows from Lemma 2.2 and Lemma 2.3 of [W].
\hfill$\Box$\vskip 5mm

\section {Proof of Theorem 1}

{\bf In this section, we assume $\Delta(\alpha,\beta) \geq6$ and the
endpoints of edges $\Gamma_P$ are with Type B labels.}

{\bf Lemma 5.1.} \ There are not two edges  which are parallel in
both $\Gamma_P$ and $\Gamma_Q$.\label{lemmaparallel}

{\bf Proof} \ The proof follows from Lemma 2.1 of
[SW].\hfill$\Box$\vskip 5mm

{\bf Lemma 5.2.} \ $\Gamma_P$ can not have $2q$ parallel edges.
\label{lemma2q}

{\bf Proof} \ Suppose $S=\{e_1, e_2, \cdots,e_{2q}\}$ is a
collection of $2q$ parallel edges joining $\partial_u P$ and
$\partial_v P$ in $\Gamma_P$, where $\partial e_i=x_i\cup y_i$.

Let $x\in\{x_1,x_2,\cdots,x_{2q},y_1,y_2,\cdots,y_{2q}\}$, give new
labels on $x$ as follows:

(1) label $x$ with $i$ if $x$ is labeled $+i$.

(2) label $x$ with $2q+1-i$ if $x$ is labeled $-i$.

These labels of $S$ give a permutation $\pi$ of $\{1, 2, \cdots,
2q\}$ defined $\pi(a) = b$ if $(a,b)$ is a label pair of an edge in
$S$. One can see that $\pi(a) = -a + s( mod 2q)$, where $s$ is a
constant. It follows that $\pi^2(a)= a$. This means if there is an
edge $e_i$ with label pair $(a,b)$, then there is a dual edge in $S$
with label pair with $(b,a)$. By Lemma 3.3, $a\neq b$. Then $S$ can
be divided into $q$ pairs, each of them consists a pair edges of
$e_k$ and $e_k^\prime$ in $S$ such that they have the same label
pair, that is they form a length 2 cycle in $\Gamma_Q$. Suppose
$e_{k_0}$ and $e_{k_0}^\prime$ is a pairs such that they form an
innermost length 2 cycle in $\Gamma_Q$. Then $e_{k_0}$ and
$e_{k_0}^\prime$ are parallel in both $\Gamma_P$ and $\Gamma_Q$,
contradicting  Lemma 5.1. \hfill$\Box$\vskip 5mm

{\bf Lemma 5.3.} \ Let $\Gamma$ be a graph embedded in a 2-sphere
with $V$($V\geq3$) vertices and $E$ edges, if $\Gamma$ contains no
1-sided disk faces and no 2-sided disk faces, then
$E\leq3V-6$.\label{lemma3v6}

{\bf Proof} \ Suppose $\Gamma$ contains $F$ faces, and all of them
have at least 3 sides, then $2E \geq 3F$. Hence $V-E+(2/3)E\geq2$
and $E\leq3V-6$. \hfill$\Box$\vskip 5mm

{\bf Lemma 5.4.} \ $p\geq5$.

{\bf Proof} \ Suppose, otherwise, that $p\leq 4$. Let
$\bar{\Gamma}_P$ be a reduced graph of $\Gamma_P$. Then
$\bar{\Gamma}_P$ has no 1-sided and no 2-sided disk-faces. Since $M$
is simple, $p>2$. By Lemma 5.3, there are at most 6 edges in
$\bar{\Gamma}_P$. Hence there is at least one vertex of
$\bar{\Gamma}_P$ which has valency at most 3. Since $\Delta \geq6$,
$\Gamma_P$ contains  $2q$ parallel edges, contradicting Lemma 5.2.
\hfill$\Box$\vskip 5mm

An $i$-collection is a collection $S=\{e_1,e_2,\cdots,e_n\}$ of
adjacent parallel  edges in $\Gamma_P$ such that each of $e_1$ and
$e_n$ has $+i$ as a signed label.

{\bf Lemma 5.5.} \ Suppose $S$ is an  $i$-collection, then the
signed labels of the endpoints of $S$ must appear as one of the
following six types.

{\bf Type I}: Each edge in $S$ has the same labels with opposite
signs. In this case, $S$ contains a virtual S-cycle labeled
$\{1,1\}$ but no S-cycle. See the following figure.

\begin{center}
\includegraphics[totalheight=4cm]{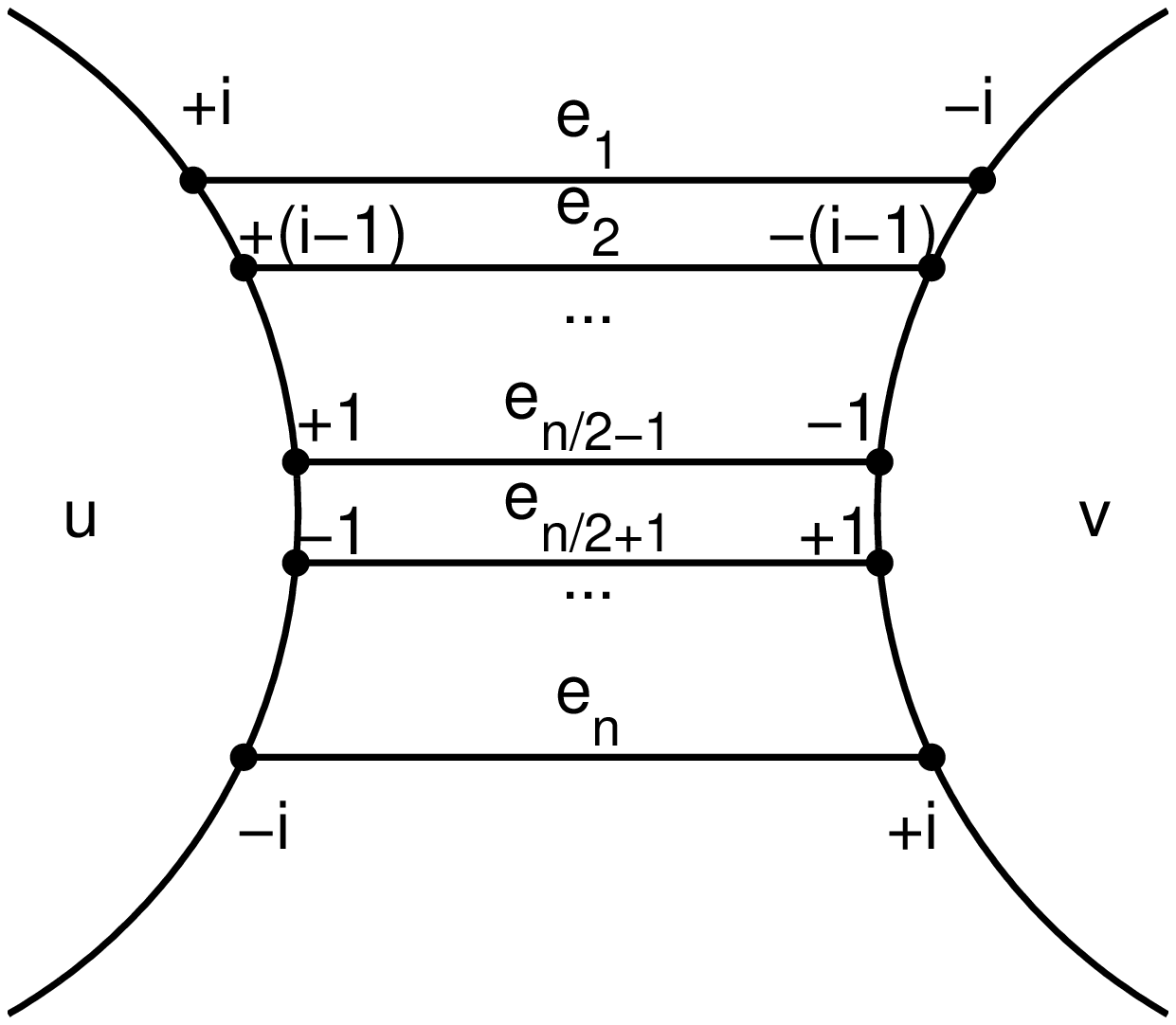}\hfill

Type I \vskip 3mm
\end{center}

{\bf Type II}: $S$ contains an S-cycle labeled $\{1,2\}$. See the
following figures.

~\hfill
\includegraphics[totalheight=4cm]{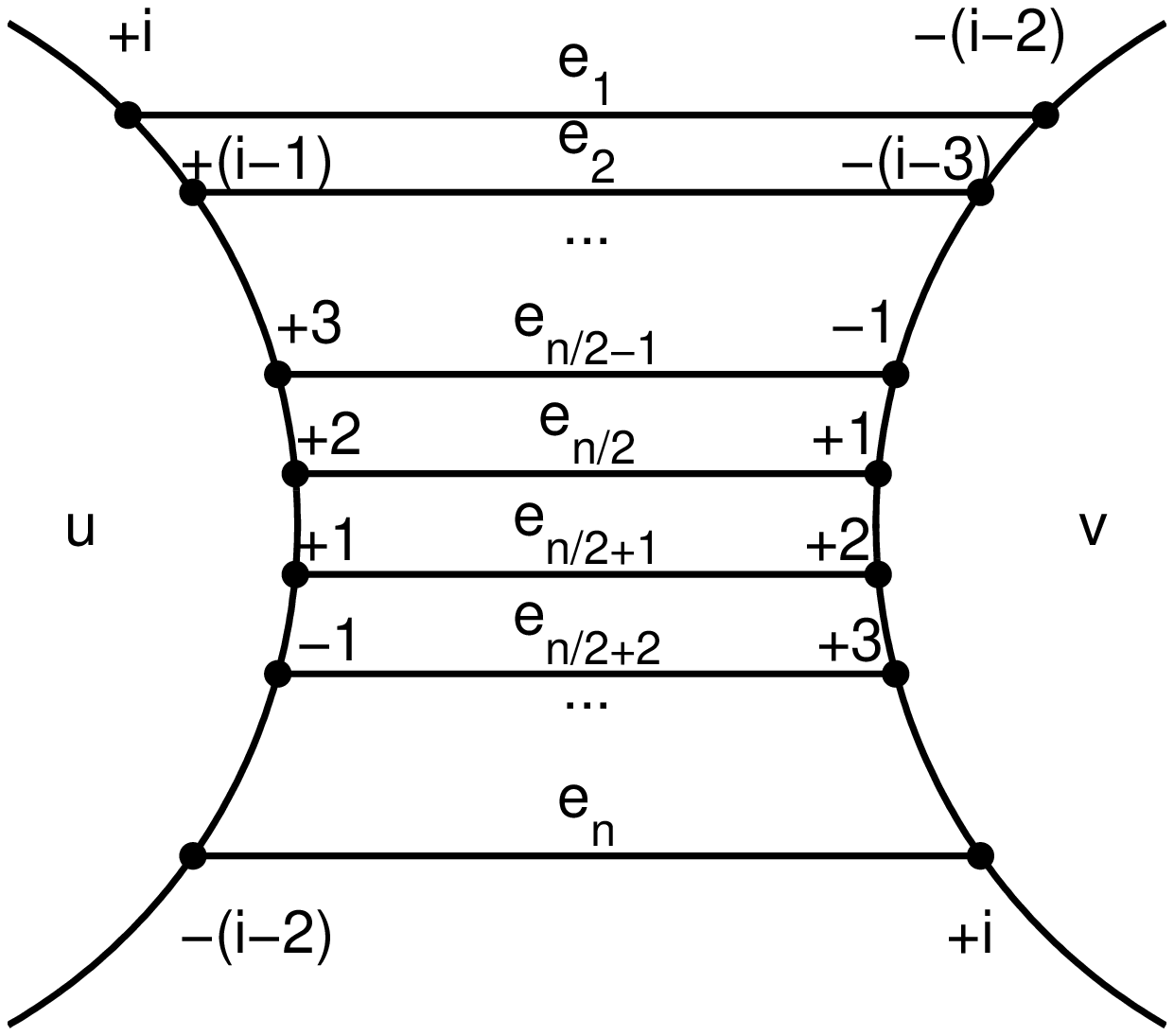}\hfill
\includegraphics[totalheight=4cm]{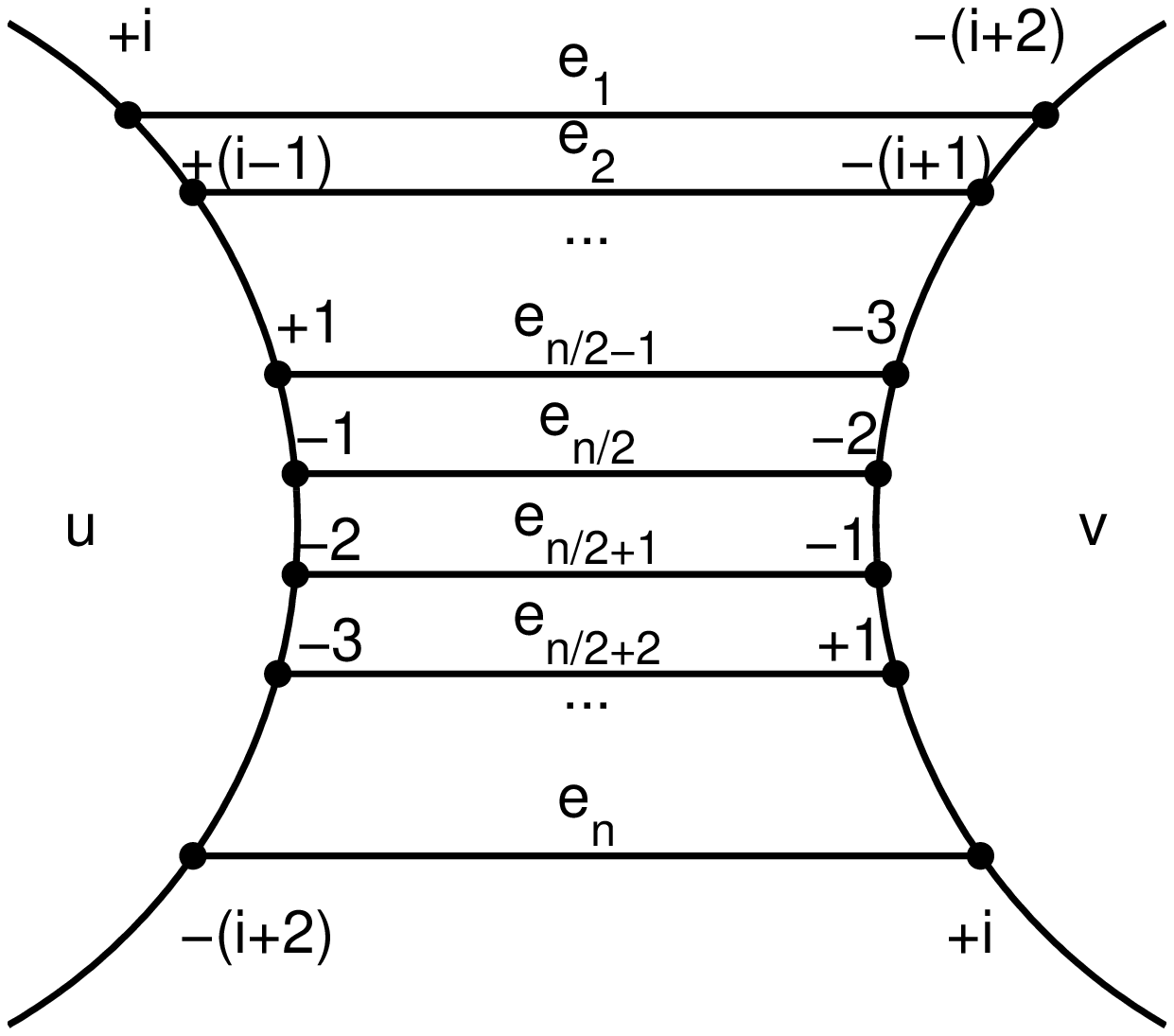}\hfill~

\begin{center} Type II(a)~~~~~~~~~~~~~~~~~~~~~~~~~~~~~~Type II(b)
\end{center}

{\bf Type III}: $n=2$, and $\{e_1, e_2\}$ is an  S-cycle labeled
$\{i,i-1\}$, where $i>1$. See the following figure.

\begin{center}
\includegraphics[totalheight=2.8cm]{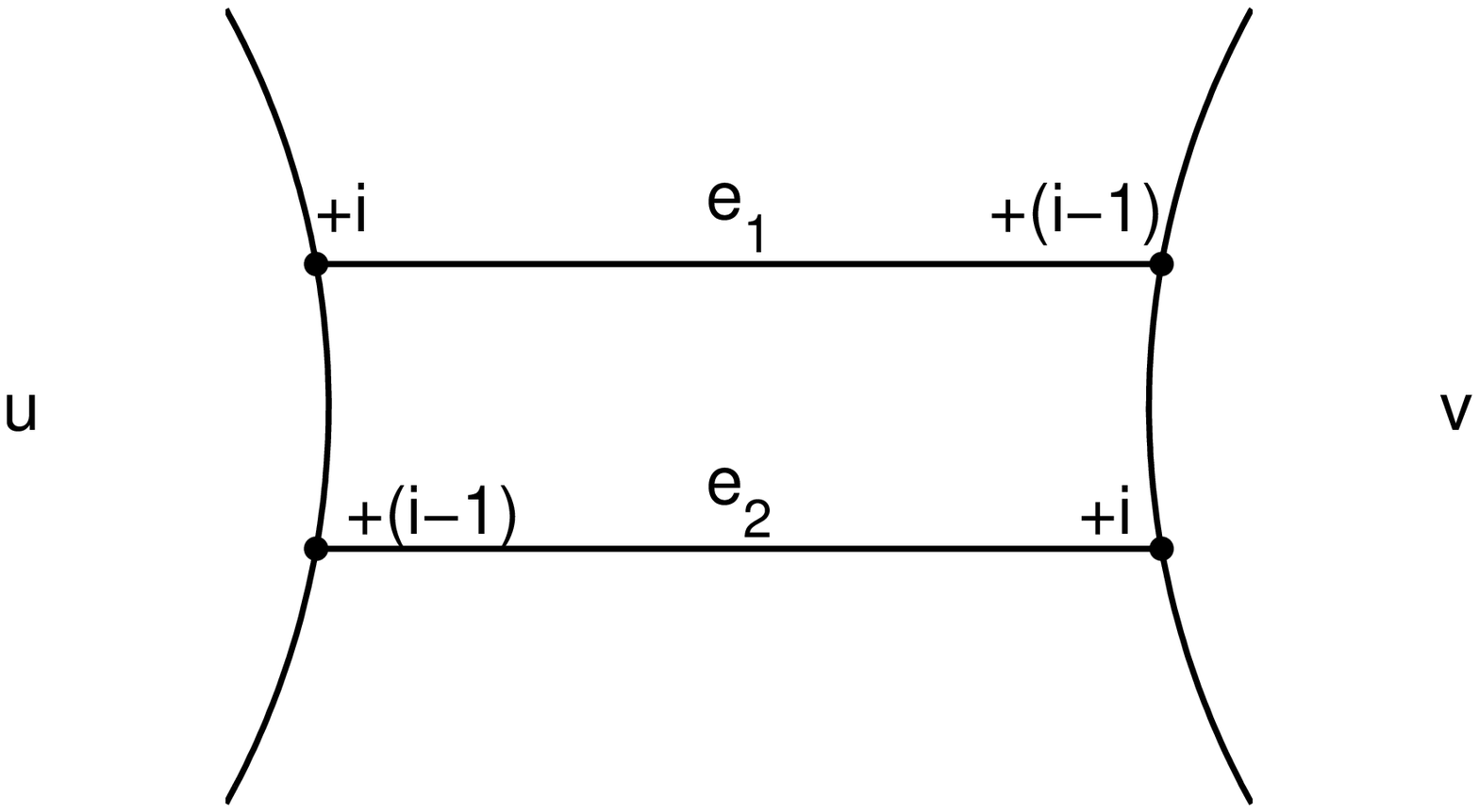}\hfill

Type III \vskip 3mm
\end{center}
{\bf Type IV}: $n=2$, and $\{e_1, e_2\}$ is an  S-cycle labeled
$\{i,i+1\}$, where $i<q$. See the following figure.
\begin{center}
\includegraphics[totalheight=2.8cm]{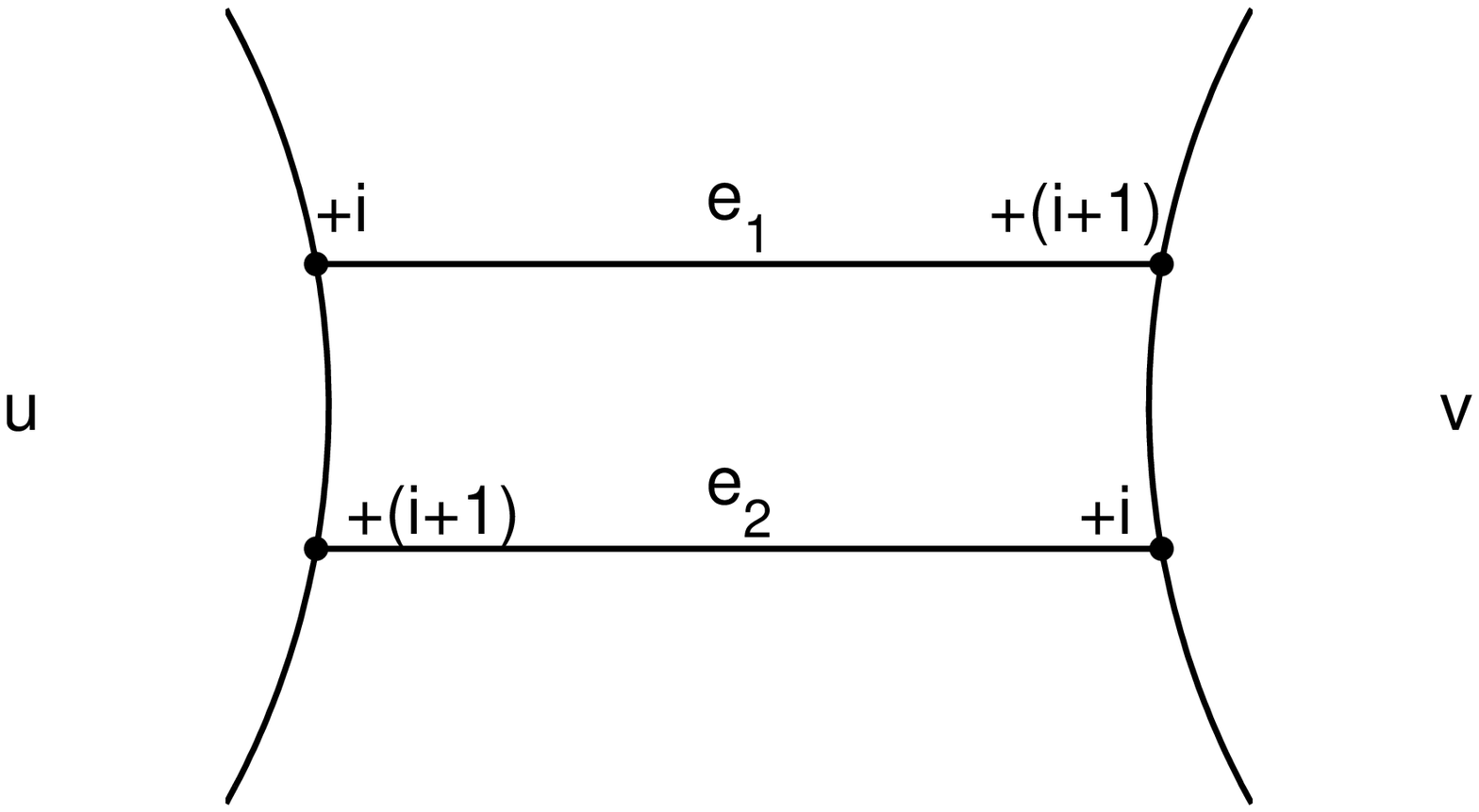}\hfill

Type IV \vskip 3mm
\end{center}

{\bf Type V}: $S$ contains an S-cycle labeled $\{q,q-1\}$. See the
following figures.

\begin{center}
\includegraphics[totalheight=4.5cm]{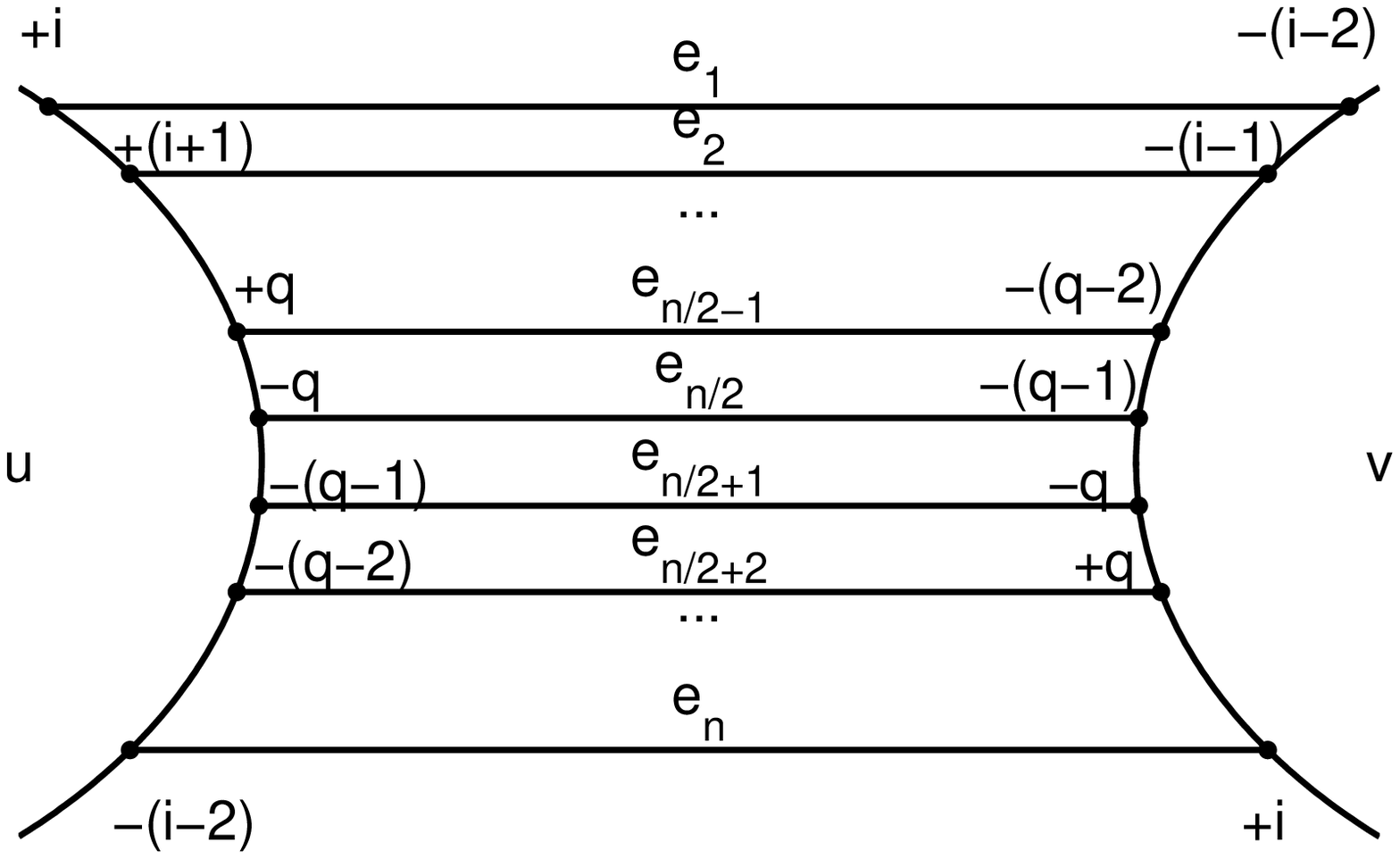}~~~~~~~~~~
\includegraphics[totalheight=4.5cm]{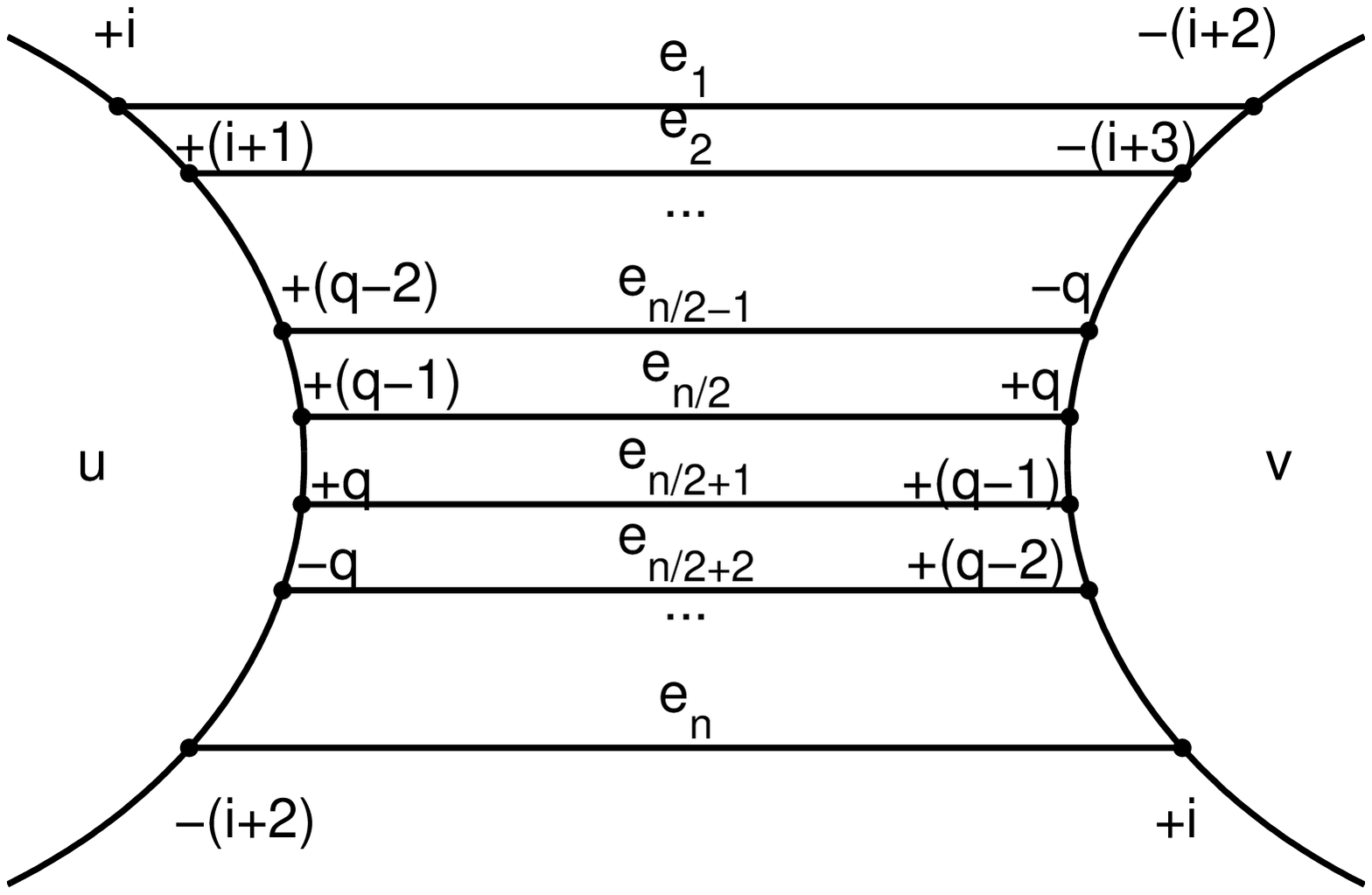}

Type V(a)~~~~~~~~~~~~~~~~~~~~~~~~~~~~~~~~~~~~~~~~~~~~Type V(b)
\end{center}

{\bf Type VI}: Each edge in $S$ has the same labels with opposite
signs. In this case, $S$ contains a virtual S-cycle labeled
$\{q,q\}$ but no S-cycle. See the following figure.

\begin{center}
\includegraphics[totalheight=4cm]{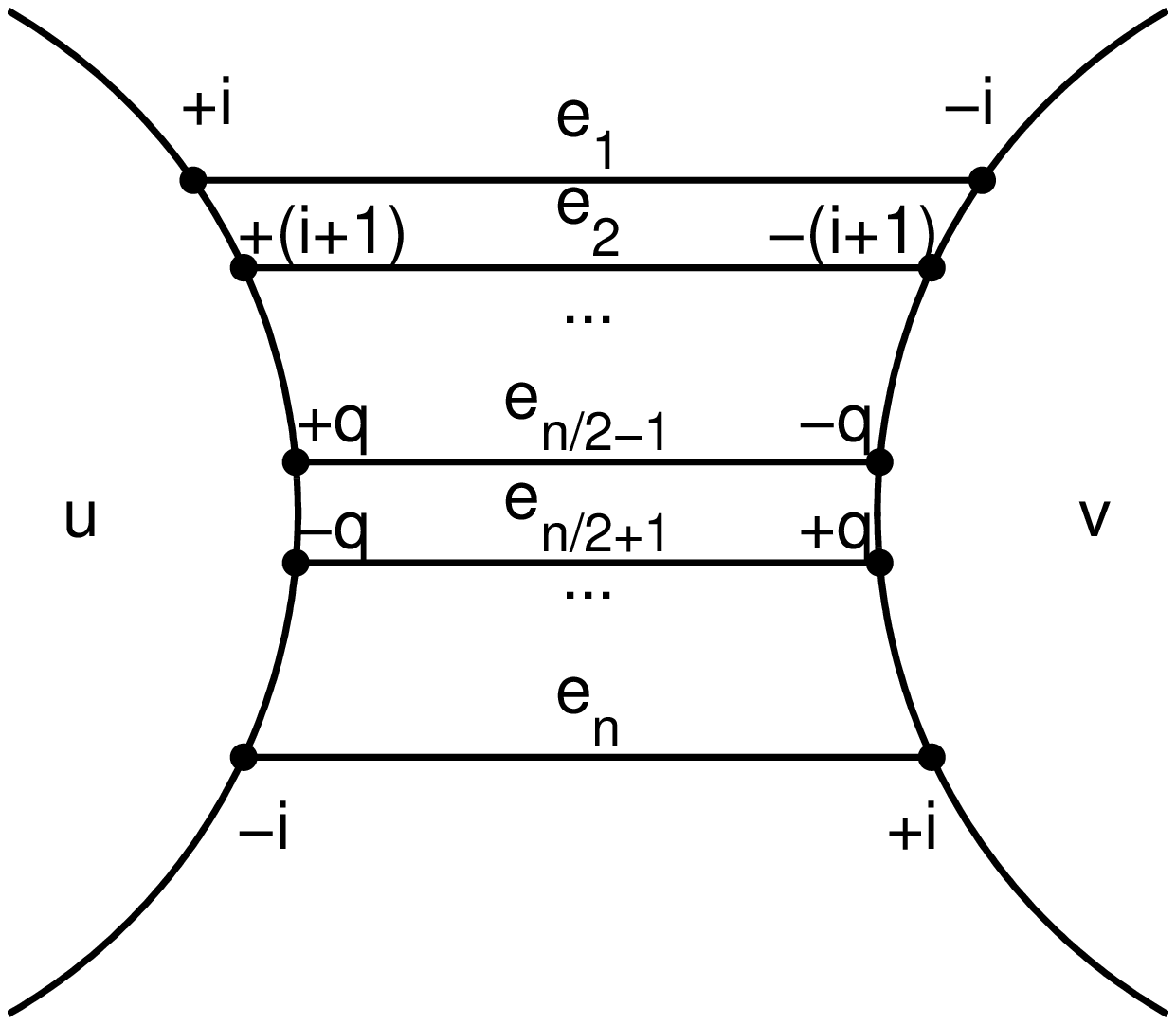}

Type VI \vskip 3mm
\end{center}

{\bf Proof} \ Assume that $\partial e_k=x_k\cup y_k$ such that
$x_k\in \partial_u P$ and $y_k \in \partial_v P$, and  $x_1$ is
labeled with $+i$. Since $S$ is an $i$-collection, by definition,
one of $x_n$ and $y_n$ is labeled with $+i$. If $x_n$ is labeled
with $+i$, then $n\geq2q$, contradicting lemma 5.2. Hence $y_n$ is
labeled with $+i$. By remark ($*$) the signed labels
$\{1,2,\cdots,q,-q,\cdots\}$ appear in the same direction in
$\Gamma_P$. Hence the signed labels of $x_{1+k}$ is the same with
the one of $y_{n-k}$ for all $k=0,1,\cdots,n$. It follows that $n$
is even; otherwise, $x_{(1+n)/2}=y_{(1+n)/2}$, contradicting Lemma
3.3.

As signed labels, we assume that $-1<+1$ and $+q<-q$.

{\bf Case 1.} \ The signed label of $x_2$ is smaller than the one of
$x_1$.

{\bf Case 1.1} \  $n=2$.

Now $S$ is a virtual S-cycle. If $x_1$ and $x_2$ are labeled with
$+1$ and $-1$, then $S$ is of type I. If $x_1$ and $x_2$ are with
$+i$ and  $+(i-1)$ for some $2\leq i\leq q$, then $S$ is of type
III.

{\bf Case 1.2 } \ $n\geq4$.

Now $\{e_{n/2-1},e_{n/2},e_{n/2+1},e_{n/2+2}\}$ is a virtual
extended S-cycle. By Lemma 4.2(2), it is not an extended S-cycle.
Hence  $\{e_{n/2},e_{n/2+1}\}$ is labeled with one of
$(1,1)$,$(1,2)$,$(q,q-1)$,$(q,q)$. Since  the signed label of $x_2$
is smaller than the one of $x_1$,  $S$ contains at least $2q$ edges
when $\{e_{n/2},e_{n/2+1}\}$ is labeled with one of $(q,q-1)$ and
$(q,q)$. Hence $S$ is one of type I and type II.

{\bf Case 2} \ the signed label of $x_2$ is bigger than the one of
$x_1$.

By the same argument as above, $S$ is one of type IV, type V and
type VI. \hfill$\Box$\vskip 5mm

{\bf The proof of Theorem 1 }

For each $1\leq i \leq q$, let $B_P^{+i}$ be a subgraph of
$\Gamma_P$ consisting all the vertices of $\Gamma_P$ and all the
edges $e$ such that one endpoint of $e$ is labeled with $+i$.

Since $\Delta(\alpha,\beta)\geq 6$, by Lemma 3.3, there are at least
$3p$ edges in $B_P^{+i}$. By Lemma 5.3, $B_P^{+i}$ contains at least
one 2-sided face. Hence there is at least one $i$-collections in
$\Gamma_P$ for each $i$. \vskip 5mm

{\bf Claim 1} For each $1\leq s\leq q-1$, $\Gamma_P$ contains no an
$s$-collection of type I and an $(s+1)$-collection of type VI
 simultaneously.\label{lemmaonesix}

{\bf Proof} \ Suppose, otherwise, that $S_1$ is an $s$-collection of
type I and $S_2$ is an $(s+1)$-collection of type VI. By the
definitions of type I and type VI, for each $i\leq s$, there are two
edges in $S_1$ with both two endpoints incident to $\partial_i Q$;
for each $j\geq s+1$, there are two edges in $S_2$ with both two
endpoints incident to $\partial_j Q$.  Hence each edge in $S_1\cup
S_2$ is a length 1 cycle in $\Gamma_Q$. This means that $\Gamma_Q$
contains a 1-sided disk-face, contradicting Lemma 2.2.
\hfill$\Box$(Claim 1)\vskip 5mm

{\bf Claim 2} \ For each $1\leq s\leq q-1$, $\Gamma_P$ contains no
an $s$-collection of type I (resp. II) and an $(s+1)$-collection of
type V(resp. VI) simultaneously.\label{lemmaonefive}

{\bf Proof} \ Suppose, otherwise, that $S_1$ is an $s$-collection of
type I and $S_2$ is an $(s+1)$-collection of type V. By the
definition of type V, all the vertices of $\partial_i Q$($i\geq
s+1$) be connected by the edges in $S_2$. By the definition of type
I, each edge in $S_1$ is a length 1 cycle which bounds two disks in
$\hat{Q}$, say $D_1$ and $D_2$. We may assume that $D_1$ is disjoint
from $\partial_i Q$ for each $i\geq s+1$. Hence $\Gamma_Q$ contains
a 1-sided disk-face, a contradiction.

Similarly, one can prove that, $\Gamma_P$ contains no an
$s$-collection of type II and an $(s+1)$-collection of type VI
simultaneously.\hfill$\Box$(Claim 2)

{\bf Claim 3} \ For each $1\leq s \leq q$, $\Gamma_P$ contains
neither $s$-collections of type II nor $s$-collections of type V.

{\bf Proof} \ Suppose, otherwise, that there is an $s$-collection of
type II for some $1\leq s\leq q$. Then $\Gamma_P$ contains an
S-cycle labeled $\{1,2\}$. By Lemma 4.3, each $i$-collection is one
of type I, type II and type VI for each $i>2$.

By Claim 1 and Claim 2, the $q$-collection is one of type I and type
II.

If a $q$-collection is of type I, then $\Gamma_P$ contains $2q$
edges, contradicting Lemma 5.2.

Assume that there is a $q$-collection of type II. Then there are two
edges connecting $\partial_1 Q$ to $\partial_2 Q$, and two edges
connected $\partial_k Q$ to $\partial_{k+2} Q$ for all $1\leq k\leq
q-2$. Each pair of the two edges as above is a length 2 cycle in
$\Gamma_Q$. Let $c$ be an innermost one of all such cycles. Then the
two edges in $c$ are parallel in both $\Gamma_P$ and $\Gamma_Q$,
contradicting Lemma 5.1.

Similarly, $\Gamma_P$ contains no $s$-collections of type
V.\hfill$\Box$(Claim 3)\vskip 5mm

{\bf Claim 4} \ For each $1\leq s \leq q$, $\Gamma_P$ contains
neither $s$-collections of III nor $s$-collections of type IV.

{\bf Proof} \ Suppose, otherwise, that there is an $s$-collection
 $S$ of type IV, then $1\leq s\leq q-1$. Note that $S$ is also an
$(s+1)$-collection of type III. By Lemma 4.3, there is no
$i$-collections of type III and VI for $i\neq s,s+1$. By Claim 3,
each $i$-collection is either of type I or type VI for $i\neq
s,s+1$. If $s=1$, then $S$ is also a $2$-collection of type II,
contradicting Claim 3. Hence $2\leq s\leq q-1$.

Since a $1$-collection of type VI contains $2q$ edges, contradicting
Lemma 5.2. Hence all the $1$-collections are of type I. By Claim 1,
all the $i$-collections are of type I for $i<s$. By Lemma 5.2, all
the $q$-collections are of type VI. By Claim 1, all the
$j$-collections are of type VI for $j>s+1$.

Suppose $S_1$ is an $(s-1)$-collection of type I, and $S_2$ is an
$(s+2)$-collection of type VI. By the definitions of type I and type
VI, there is  a length 1 cycle incident to $\partial_i Q$ for each
$i\neq s,s+1$. Since two edges in $S$ connect $\partial_s Q$ to
$\partial_{s+1} Q$, So $\Gamma_Q$ contains a 1-sided face in
$\Gamma_Q$, it is a contradiction to Lemma 2.2. \hfill$\Box$(Claim
4)

By the Claim 3 and Claim 4. all the $s$-collections are of type I or
type VI for $1\leq s\leq q$. since all $1$-collections are of type
I, and all $q$-collections are of type VI. Then we can find some $k$
such that there are a $k$-collection of type I and a
$k+1$-collection of type VI, contradicting Claim 1. \hfill$\Box$
(Theorem 1)

\noindent Mingxing Zhang \\
Department of Applied Mathematica, Dalian University of Technology,
\\Dalian, China, 116022 \\Email: star.michael@263.net

\vskip 5mm
\noindent Ruifeng Qiu\\
Department of Applied Mathematica, Dalian University of
Technology,\\
Dalian, China, 116022  \\Email: qiurf@dlut.edu.cn

\vskip 5mm
\noindent Yannan Li \\
Department of Applied Mathematica, Dalian University of
Technology,\\ Dalian, China, 116022\\ Email: yn\_lee79@yahoo.com.cn

\end{document}